\documentclass[twoside,a4paper,reqno,11pt]{amsart} 
\usepackage{amsfonts, amsbsy, amsmath, amssymb, latexsym,hyperref, mathtools, bold-extra, url}
\usepackage{mathrsfs,array}
\usepackage[top=24mm,right=28mm,bottom=24mm,left=28mm]{geometry}
\usepackage{stmaryrd}
\usepackage{bm}
\usepackage[pdftex]{color,graphicx}

\newcommand{\e}{\epsilon}
 
\renewcommand{\l}{\lambda} \renewcommand{\O}{\Omega}

 \renewcommand{\to}{\rightarrow}
 \newcommand{\s}{\sigma}
 
 \newcommand{\C}{\mathcal{C}}

\newcommand{\la}{\langle}
\newcommand{\ra}{\rangle}

\newcommand{\leqs}{\leqslant}
\newcommand{\geqs}{\geqslant}

 \newcommand{\vs}{\vspace{3mm}}

\makeatletter
\newcommand{\imod}[1]{\allowbreak\mkern4mu({\operator@font mod}\,\,#1)}
\makeatother

\newtheorem{theorem}{Theorem} 
\newtheorem*{theorem*}{Theorem} 
\newtheorem*{conj*}{Conjecture}

\newtheorem{thm}{Theorem}[section] 
\newtheorem{prop}[thm]{Proposition}

\theoremstyle{definition}
\newtheorem{rem}[thm]{Remark}
\newtheorem{remk}{Remark}

\newtheorem{defn}[thm]{Definition}

\begin{document}

\author{Timothy C. Burness}
\address{T.C. Burness, School of Mathematics, University of Bristol, Bristol BS8 1UG, UK}
\email{t.burness@bristol.ac.uk}

\author{Lei Wang}
\address{L. Wang, School of Mathematics and Statistics, Yunnan University, Kunming 650091, Yunnan, P. R. China}
\address{Yunnan Key Laboratory of Applied Mathematics, Yunnan University, Kunming 650091, Yunnan, P. R. China}
\email{wanglei@ynu.edu.cn}

\title[Orders of irreducible subgroups]{On the orders of irreducible subgroups \\ of finite classical groups}

\begin{abstract}
Let $G$ be a finite almost simple classical group over $\mathbb{F}_q$ with socle $G_0$ and natural module $V$ of dimension $n$. Let $H$ be a core-free maximal subgroup of $G$. By a well-known theorem of Aschbacher, either $H$ is contained in one of eight subgroup collections that are defined in terms of the geometry of $V$, or $H$ is almost simple and the embedding of $H$ in $G$ arises from an absolutely irreducible representation of a covering group of the socle of $H$. We write $\mathcal{S}$ for the latter collection of irreducibly embedded almost simple subgroups. Determining the subgroups in $\mathcal{S}$ is an intractable problem, in general, which has only been fully resolved in low dimensions. However, many applications  require some control on the subgroups in $\mathcal{S}$, expressed  in terms of their order, or properties of their elements, for example. Extending a theorem of Liebeck from the 1980s, and building on later work by Burness, Guralnick and Saxl, in this paper we determine all the subgroups $H \in \mathcal{S}$ satisfying the lower bound $|H| \geqs  q^n$. We also highlight a recent application of our main result to a problem in finite geometry. 
\end{abstract}

\date{\today}

\maketitle

\section{Introduction}\label{s:intro}

Let $G$ be a finite almost simple classical group over $\mathbb{F}_q$ with socle $G_0$ and let $V$ be the natural module for $G_0$. Write $q = p^f$ and $n = \dim V$, where $p$ is a prime.

The main theorem on the subgroup structure of finite classical groups is due to Aschbacher. In \cite{asch}, eight collections of subgroups of $G$ are introduced, denoted by $\C_1, \ldots, \C_8$ and defined in terms of the geometry of $V$. For instance, the collection $\C_1$ comprises the stabilisers in $G$ of suitable subspaces (or pairs of subspaces) of $V$. And the subgroups in $\C_2$ are the stabilisers of certain direct sum decompositions of $V$, and so on. Aschbacher's main theorem then states that if $H$ is any core-free maximal subgroup of $G$, then either $H$ is contained in one of the $\C_i$ collections, or $H$ is almost simple and the embedding of $H$ in $G$ is afforded by an absolutely irreducible representation $\rho:L \to {\rm GL}(V)$, where $L$ is a quasisimple covering group of the socle of $H$. In the latter situation, we say that $H$ is contained in the collection $\mathcal{S}$ and we note that the image $\rho(L)$ must satisfy several additional conditions that are designed to ensure that $H$ is not contained in one of the geometric $\C_i$ collections (see Definition \ref{d:sdef}). 

In \cite{KL}, Kleidman and Liebeck present detailed information on the geometric subgroups comprising each $\C_i$ collection. In particular, they calculate the order of every subgroup of this form and more importantly, assuming $n \geqs 13$, they determine all of the maximal subgroups in $\C_i$ up to conjugacy in $G$. This is complemented by the more recent work of Bray, Holt and Roney-Dougal \cite{BHR}, which completely determines the maximal subgroups, up to conjugacy, in the low-dimensional classical groups with $n \leqs 12$. The latter work in \cite{BHR} also deals with the subgroups in $\mathcal{S}$, and the analysis has subsequently been extended to the groups with $n \in \{13,14,15\}$ by Schr\"{o}der \cite{Schroder} and $n \in \{16,17\}$ by Rogers \cite{Rogers}. However, the problem of determining the maximal subgroups in $\mathcal{S}$ (even up to isomorphism) for all $n$ is intractable.

The latter observation has an impact on a wide range of applications. For instance, by combining the O'Nan-Scott Theorem with the Classification of Finite Simple Groups, we obtain a powerful strategy that can sometimes be used to reduce a general problem on finite primitive permutation groups to the case where $G \leqs {\rm Sym}(\O)$ is an almost simple classical group with point stabiliser $H$. And since $H$ is a core-free maximal subgroup of $G$, we may be able to use  Aschbacher's theorem to reduce further to the situation where $H$ is contained in $\mathcal{S}$. At this point, additional  information concerning certain properties of $H$ will usually be needed in order to reach the desired conclusion, despite the fact that we do not know the subgroups comprising $\mathcal{S}$, in general.

For some applications, it is sufficient to have restrictions on the order of $H$, even if we cannot determine its precise structure. For example, in the permutation group setting, upper bounds on \emph{fixed point ratios} play a central role in the application of probabilistic methods to a wide range of problems, from group generation to the study of base sizes (we refer the reader to the survey article \cite{Bur18} for more details). Given an element $x \in G$, we write
\[
{\rm fpr}(x,\O) = \frac{|x^G \cap H|}{|x^G|}
\]
for the fixed point ratio of $x$, which is simply the proportion of points in $\O$ fixed by $x$. So in order to derive an upper bound on ${\rm fpr}(x,\O)$, we need to combine an upper bound on $|x^G \cap H|$ with a lower bound on $|x^G|$. Of course, we can obtain a crude upper bound on the numerator by taking $|x^G \cap H| \leqs |H|$, so long as we have some control on the order of $H$. 

In the special case where $G$ is a classical group as above and $H$ is a subgroup in the collection $\mathcal{S}$, it is natural to seek an upper bound on $|H|$ in terms of the defining parameters $n$ and $q$ (for example, this can then be combined with a lower bound on $|x^G|$, which can also be expressed in terms of $n$ and $q$). In practice, it may be desirable to derive a more refined statement, which states that the pair $(G,H)$ either belongs to a known list of candidate cases (and for each one, we can work directly with the embedding of $H$ in $G$), or $|H| < f(n,q)$ for some explicit function in $n$ and $q$.

The first result in this direction is due to Liebeck \cite{L85} and it has been very useful in a wide range of applications. In \cite[Theorem 4.2]{L85}, all the pairs $(G,H)$ such that $H \in \mathcal{S}$ and $|H| \geqs f(n,q)$ are determined, where 
\[
f(n,q) = \left\{\begin{array}{ll}
q^{4n+8} & \mbox{if $G_0 = {\rm U}_n(q)$ is a unitary group,} \\
q^{2n+4} & \mbox{otherwise.}
\end{array}\right.
\]

For certain applications, it is desirable to extend Liebeck's theorem by weakening the lower bound on $|H|$ at the expense of expanding  the list of possibilities for $(G,H)$. For example, the main theorem of  \cite{BGS} determines the exact base size of every finite almost simple primitive permutation group with classical socle and point stabiliser in $\mathcal{S}$, and the probabilistic proof relies heavily on fixed point ratio estimates. In order to simplify the proof for unitary groups, Liebeck's lower bound of $|H| \geqs q^{4n+8}$ in this case was extended to $|H| \geqs q^{2n+4}$ via \cite[Theorem 2.10]{BGS}, bringing it in line with the existing lower bound for all the other classical groups.

With a view towards further applications, the main aim of this paper is to significantly extend the earlier work on this problem by replacing the lower bound $|H| \geqs q^{2n+4}$ in \cite[Theorem 2.10]{BGS} by $|H| \geqs q^n$.

\begin{theorem}\label{t:main}
Let $G$ be a finite almost simple classical group over $\mathbb{F}_q$ with socle $G_0$ and let $n$ be the dimension of the natural module $V$ for $G_0$. Let $H \in \mathcal{S}$ be a maximal subgroup of $G$ and let $H_0$ be the socle of $H \cap G_0$. Then one of the following holds:
\begin{itemize}\addtolength{\itemsep}{0.2\baselineskip}
\item[{\rm (i)}] $H_0 = A_d$, $q=p$ is a prime and $V$ is the fully deleted permutation module for $H_0$.
\item[{\rm (ii)}] $(G_0,H_0)$ is one of the cases recorded in Table \ref{tab:nondef1} or \ref{tab:def1}.
\item[{\rm (iii)}] $|H| < q^n$.
\end{itemize}
\end{theorem}

\begin{remk}\label{r:main}
Let us record some comments on the statement of Theorem \ref{t:main}.
\begin{itemize}\addtolength{\itemsep}{0.3\baselineskip}
\item[{\rm (a)}] In the statement, and throughout the paper, we assume $G_0$ is one of the following:
\[
\hspace{12mm} {\rm L}_n(q) \, (n \geqs 2), \; {\rm U}_n(q) \, (n \geqs 3), \; {\rm PSp}_n(q)' \, (n \geqs 4), \; {\rm P\O}_n^{\e}(q) \, (n \geqs 7),
\]
where we adopt the standard notation for simple groups from \cite{KL} (in addition, we will sometimes write ${\rm L}_n(q) = {\rm L}_n^{+}(q)$ and ${\rm U}_n(q) = {\rm L}_n^{-}(q)$). This is a valid assumption due to the existence of isomorphisms between some of the low-dimensional classical groups (see \cite[Proposition 2.9.1]{KL}). 

Similarly, we assume that $G_0 \ne {\rm L}_2(5), {\rm L}_2(7), {\rm Sp}_4(2)', {\rm PSp}_4(3)$ since
\[
\hspace{12mm}  {\rm L}_2(5) \cong {\rm L}_2(4), \; {\rm L}_2(7) \cong {\rm L}_3(2),\; {\rm Sp}_4(2)' \cong {\rm L}_2(9), \; {\rm PSp}_4(3) \cong {\rm U}_4(2).
\]

\item[{\rm (b)}] In (i), $H_0 = A_d$ is an alternating group and $V$ is the fully deleted permutation module for $H_0$ over $\mathbb{F}_p$. We refer the reader to \cite[Section 5.3]{KL} for details on the  construction of this familiar module, together with the list of pairs $(G,H)$ that arise (also see \cite[Table 3]{BGS}). For example, $G = {\rm Sp}_{n}(2)$ and $H = S_{n+2}$ with $n \equiv 0 \imod{4}$ is one possibility, in which case there is clearly no upper bound of the form $|H| < q^{f(n)}$ (valid for all $n$) for any linear function $f(n) = an+b$. This explains why the cases involving the fully deleted permutation module need to be recorded separately in Theorem \ref{t:main}.

\item[{\rm (c)}] In (ii), we divide the relevant pairs $(G_0,H_0)$ into two parts, according to whether or not $H_0 \in {\rm Lie}(p)$, where ${\rm Lie}(p)$ is the set of finite simple groups of Lie type in characteristic $p$. The pairs with $H_0 \not\in {\rm Lie}(p)$ are recorded in Table \ref{tab:nondef1}, while the remaining pairs  appear in Table \ref{tab:def1}.

\item[{\rm (d)}] In Table \ref{tab:nondef1} for the case $(G_0,H_0) = ({\rm L}_4(13),{\rm U}_4(2))$ we have $G = G_0.\la \gamma \ra$, where $\gamma$ is an involutory graph automorphism with $C_{G_0}(\gamma) = {\rm PSp}_4(13)$. And in the final row with $(G_0,H_0) = ({\rm L}_2(9),A_5)$ we write $\phi$ for an involutory field automorphism of $G_0$. 

\item[{\rm (e)}] Let $H_0 \in {\rm Lie}(p)$ be one of the cases in Table \ref{tab:def1} and recall that we may view $V$ as an absolutely irreducible module for a quasisimple covering group $L$ of $H_0$. If $k$ is the algebraic closure of $\mathbb{F}_p$, then $\bar{V} = V \otimes k$ is an irreducible module for the ambient simple algebraic group $X$ corresponding to  $L$ and we write $\l$ for the highest weight of $\bar{V}$ as a $kX$-module, expressed in terms of a set of fundamental dominant weights $\{\l_1, \ldots, \l_r\}$ for $X$, which we label in the usual manner (see \cite{Bou}). In the third column of Table \ref{tab:def1}, we record the highest weight $\l$ of $\bar{V}$, up to quasiequivalence.

\item[{\rm (f)}] In Table \ref{tab:def1} we adopt the following notation:

\vspace{1mm}

\begin{itemize}\addtolength{\itemsep}{0.2\baselineskip}
\item In the first row, we set $z=1$ if $\e=+$ and $z=d-1$ if $\e=-$. 

\item For $G_0 = {\rm P\O}_{n}^{\e}(q)$ and $H_0 = {\rm L}_d^{\pm}(q)$ with $\l = \l_1+\l_{d-1}$, we set $c_1 = 2$ if $p$ divides $d$, otherwise $c_1=1$. 

\item For $G_0 = {\rm P\O}_{n}^{\e}(q)$ and $H_0 = {\rm PSp}_d(q)$ with $\l = \l_2$, we set $c_2 = 2$ if $p$ divides $d/2$, otherwise $c_2=1$.
\end{itemize}

\vspace{1mm}

\item[{\rm (g)}] For the groups with $n \geqs 18$, we are not claiming that every pair $(G_0,H_0)$ appearing in Tables \ref{tab:nondef1} and \ref{tab:def1} corresponds to a maximal subgroup $H \in \mathcal{S}$. In every case, there is an irreducibly embedded subgroup of the given form, but additional conditions may be required to ensure maximality. By appealing to the relevant tables in \cite{BHR,Rogers,Schroder}, we have omitted any candidate pairs $(G_0,H_0)$ that do not give rise to a maximal subgroup $H \in \mathcal{S}$ when $n \leqs 17$. And some candidate pairs $(G_0,H_0)$ with $n \geqs 18$ have also been excluded because $H$ is clearly non-maximal. For example, if $G_0 = \O_{32}^{+}(q)$, $p=2$ and $H_0 = {\rm Sp}_{10}(q)$ with $\l = \l_5$, then $H_0 < \O_{12}^{+}(q) < G_0$ and $H$ is non-maximal.

\item[{\rm (h)}] In Table \ref{tab:def1}, there are several cases where $G_0 = {\rm P\O}_n^{\e}(q)$ with $n$ even and the sign $\e = \pm$ has not been specified. In general, there is no straightforward or uniform way to determine the sign, which will typically depend on various conditions involving the defining parameters for $H_0$. For example, we refer the reader to Remark \ref{r:sl1} for the special case where $H_0 = {\rm L}_d(q)$ and $\l = \l_{1}+\l_{d-1}$.
\end{itemize}
\end{remk}

{\small
\begin{table}
\[
\begin{array}{lll} \hline
G_0 & H_0 & \mbox{Conditions} \\ \hline
{\rm U}_{40}(q) & {\rm PSp}_8(3) & q = 2\\
{\rm P\O}_{28}^{-}(q) & {\rm Ru}  & q = 2\\
{\rm U}_{27}(q) & \O_7(3),\, {\rm Fi}_{22}  & q = 2 \\
{\rm L}_{24}(q) & {\rm PSp}_4(7)  & q = 2 \\
{\rm P\O}_{24}^{+}(q) & {\rm Co}_1  & q = 2,3,5 \\
\O_{23}(q) & {\rm Co}_2 &  q = 3 \\
{\rm PSp}_{22}(q) & {\rm Co}_3  & q = 2 \\
{\rm P\O}_{22}^{+}(q) & {\rm Co}_3 &  q = 3\\
 & {\rm Co}_2,\, {\rm McL} &  q  = 2 \\
\O_{21}(q) & {\rm U}_6(2) & q = 3 \\
{\rm PSp}_{20}(q) & {\rm HS}  & q = 2 \\
{\rm P\O}_{20}^{+}(q) & {\rm U}_4(3)  & q = 2\\
{\rm L}_{16}^{\e}(q) & A_{12}  & (q,\e) = (3,+) \\
 & A_{11}  & (q,\e) = (2,-) \\
{\rm U}_{14}(q) & {\rm PSp}_6(3)  & \mbox{$q=5$, $G = G_0.2$} \\
{\rm P\O}_{14}^{-}(q) & G_2(3) & q = 2 \\
{\rm U}_{13}(q) & {\rm PSp}_{6}(3)  &  q = 2,5 \\
{\rm L}_{12}^{\e}(q) & {\rm Suz} & (q,\e) = (2,-),(5,-),(7,+) \\
{\rm PSp}_{12}(q) & G_2(4)  & q = 5 \\
 & {\rm Suz}  & q = 3 \\
 & {\rm L}_2(25) & q = 2 \\
{\rm P\O}_{12}^{-}(q) & {\rm L}_3(3)  & q = 2 \\
{\rm L}_{11}(2) & {\rm L}_{2}(23),\, {\rm M}_{24}  & q = 2 \\
{\rm L}_{10}^{\e}(q) & {\rm M}_{12}  & (q,\e) = (3,+) \\
 & {\rm M}_{22}  & \mbox{$(q,\e) = (2,+), (3,-)$} \\
{\rm PSp}_{10}(q) & {\rm U}_5(2) & q = 3 \\
{\rm P\O}_{10}^{-}(q) & {\rm M}_{12}  & q = 2 \\
{\rm U}_9(q) & {\rm L}_2(19),\, {\rm J}_3  & q = 2 \\
{\rm PSp}_8(q) & {\rm L}_2(17) & q = 2 \\
{\rm P\O}_{8}^{+}(q) & \O_8^{+}(2)  & q = 3,5,7 \\
& & \mbox{$q = 11$, $G \ne G_0$} \\
&& \mbox{$q = 13$, $G = G_0.S_3$} \\
\O_7(q) & {\rm Sp}_6(2)  & q =3,5,7 \\
{\rm L}_6^{\e}(q) & {\rm U}_4(3)  & (q,\e) = (2,-), (5,-),(7,+),(11,-),(13,+) \\
 & {\rm M}_{12}  & (q,\e) = (3,+) \\
 & {\rm L}_3(4)  & (q,\e) = (5,-) \\
 & {\rm M}_{22}  & (q,\e) = (2,-) \\
{\rm PSp}_6(q) & {\rm J}_2  & q  = 5,9 \\
 & {\rm L}_2(13)  & q = 3 \\
& {\rm U}_3(3)  & q = 2 \\
{\rm L}_5^{\e}(q) & {\rm L}_2(11) & (q,\e) = (2,-),(3,+) \\
& {\rm U}_4(2) & (q,\e) = (5,-),(7,+) \\
& {\rm M}_{11} & (q,\e)= (3,+) \\
{\rm L}_4^{\e}(q) & A_7 & (q,\e) = (2,+),(3,-),(5,-) \\
& {\rm U}_4(2) & (q,\e) = (5,-),(7,+),(11,-) \\
& & \mbox{$(q,\e) = (13,+)$, $G = G_0.\la \gamma \ra =  G_0.2$} \\
& {\rm L}_3(4) & q = 3 \\
{\rm PSp}_4(q) & A_6 & \mbox{$q = 5$, $G = G_0.2$}  \\
& A_7 & q = 7 \\
{\rm L}_3^{\e}(q) & A_6 & (q,\e) = (4,+),(5,-) \\
& A_7 & (q,\e) = (5,-) \\
& {\rm L}_2(7) & (q,\e) = (3,-),(5,-) \\
{\rm L}_2(q) & A_5 & \mbox{$q = 9$, $G = G_0.\la \phi \ra = G_0.2$}  \\ \hline
\end{array}
\]
\caption{The cases $(G_0,H_0)$ arising in Theorem  \ref{t:main}(ii) with $H_0 \not\in {\rm Lie}(p)$}
\label{tab:nondef1}
\end{table}}

{\small
\begin{table}
\[
\begin{array}{lllll} \hline
G_0 & H_0 & \l & n & \mbox{Conditions} \\ \hline
{\rm L}_n^{\e}(q) & {\rm L}_d(q^2) & \l_1 + q\l_{z} & d^2  & \mbox{$d \geqs 3$, $\e=\pm$} \\
& {\rm L}_d^{\e}(q) & \l_2 & d(d-1)/2 &  d \geqs 5 \\
& & 2\l_1 & d(d+1)/2 &  d \geqs 3,\, p \geqs 3 \\
& & \l_3 & d(d-1)(d-2)/6  & d \in \{7,8\} \\
& {\rm P\O}_d^{\e}(q) & \l_{d/2-1} & 2^{d/2-1}  & d \in \{10,14\} \\
& E_6^{\e}(q) & \l_1 & 27 & \\
& & & &\\
{\rm PSp}_n(q) & {\rm L}_d^{\pm}(q) & \l_1+\l_{d-1} & d^2-2  & \mbox{$d \geqs 6$, $d \equiv 2 \imod{4}$, $p = 2$} \\
& & \l_3 & 20 & \mbox{$d = 6$, $p \geqs 3$} \\
& {\rm L}_d(q^3) & (1+q+q^2)\l_1 & 8 &  d = 2, \, p \geqs 3 \\
& {\rm PSp}_d(q) & \l_2 & d(d-1)/2-2 &  d \geqs 12,\, d \equiv 4 \imod{8}, \, p = 2 \\
& & \l_3 & 14 &  d = 6, \, p \geqs 3 \\
& {\rm P\O}_d^{+}(q) & \l_{d/2-1} & 2^{d/2-1} &  d = 12, \, p \geqs 3 \\
& \O_d(q) & \l_{(d-1)/2} & 2^{(d-1)/2} &  d = 13, \, p \geqs 3 \\ 
& E_7(q) & \l_7 & 56 &  p \geqs 3 \\
& & \l_1 & 132 & p = 2 \\
& G_2(q) & \l_1 & 6 & p = 2 \\ 
& {}^2B_2(q) & \l_1 & 4 & \mbox{$q = 2^f$, $f \geqs 3$ odd} \\
& & & & \\
{\rm P\O}_n^{\e}(q) & {\rm L}_d^{\pm}(q) & \l_1+\l_{d-1} & d^2-c_1 & \mbox{$d \geqs 3$, $d\not \equiv 2 \imod{4}$ or $p\geqs 3$} \\
& & \l_3 & 20 & \mbox{$d = 6$, $p = 2$, $\e=+$} \\ 
& {\rm PSp}_d(q^2) & (1+q)\l_1 & d^2 & \mbox{$d \geqs 4$, $(d,p) \ne (4,2)$} \\
& {\rm PSp}_d(q) & \l_2 & d(d-1)/2-c_2 & d \geqs 6,\, \mbox{$d \not\equiv 4 \imod{8}$ or $p \geqs 3$} \\
& & 2\l_1 & d(d+1)/2 & \mbox{$d \geqs 4$, $p \geqs 3$} \\
& & \l_{d/2} & 2^{d/2} & d \in \{ 6, 8, 12\}, \, p = 2, \, \e=+ \\ 
& {\rm P\O}_d^{+}(q) & \l_{d/2-1} & 2^{d/2-1} & d = 12, \, p = 2,\, \e = + \\
& {\rm P\O}_d^{-}(q_0) & \l_{d/2-1} & 2^{d/2-1} & d = 8, \, q = q_0^2,\, \e=+ \\
& {\rm P\O}_d^{\pm}(q) & \l_2 & d(d-1)/2 & \mbox{$d \geqs 8$, $p \geqs 3$} \\
& \O_d(q) & \l_2 & d(d-1)/2 & \mbox{$d \geqs 7$, $p \geqs 3$, $|G:G_0|>2$}\\
& & \l_{(d-1)/2} & 2^{(d-1)/2} & d \in \{7,9\}, p \geqs 3,\, \e = + \\ 
& E_8(q) & \l_8 & 248 & \\
& E_7(q) & \l_7 & 56 & p = 2,\, \e = + \\
& & \l_1 & 133 & p \geqs 3 \\
& E_6^{\delta}(q) & \l_2 & 78-\delta_{3,p} & \\ 
& F_4(q) & \l_4 & 26-\delta_{3,p} &  \\
& & \l_1 & 26(2-\delta_{2,p}) & \\
& {}^3D_4(q) & \l_2 & 28 & p \geqs 3 \\
& {}^3D_4(q_0) & \l_1 & 8 & q=q_0^3, \, \e = + \\
& G_2(q) & \l_1 & 7 & p \geqs 3 \\
& & \l_2 & 7(2-\delta_{3,p}) & \mbox{$p \geqs 3$, $\e=+$ if $q \equiv 1,7 \imod{12}$,} \\
& & & & \mbox{$\e = -$ if $q \equiv 5,11 \imod{12}$} \\ 
\hline
\end{array}
\]
\caption{The cases $(G_0,H_0)$ arising in Theorem  \ref{t:main}(ii) with $H_0 \in {\rm Lie}(p)$}
\label{tab:def1}
\end{table}}

Our proof of Theorem \ref{t:main} naturally divides into two parts, according to whether or not $H_0$ is in ${\rm Lie}(p)$. And by appealing to the detailed information on maximal subgroups in \cite{BHR,Rogers, Schroder}, it is a straightforward task to read off all the pairs $(G_0,H_0)$ with $|H| \geqs q^n$ when $n \leqs 17$. So most of our work involves the high-dimensional groups with $n \geqs 18$, where we do not know all the maximal subgroups in $\mathcal{S}$, in general.

The cases with $H_0 \not\in {\rm Lie}(p)$ are handled in Section \ref{s:notlie} and we consider each possibility for $H_0$ in turn. If $H_0$ is a sporadic simple group, then by considering the order of the Monster group we immediately reduce the problem to the groups with $n < 180$, which then allows us to appeal to work of Hiss and Malle. In \cite{HM,HM2}, all the absolutely irreducible representations of dimension at most $250$ of quasisimple finite groups are determined, excluding groups of Lie type in defining characteristic. In addition, we are also provided with information on the associated Frobenius-Schur indicators and Brauer character fields, which we can use to determine $G_0$. 

We adopt a similar approach when $H_0 \not\in {\rm Lie}(p)$ is an alternating group or a simple group of Lie type. For the alternating groups, we apply results of James \cite{James83} and Wagner \cite{Wagner77} to reduce to the groups with $n<75$, which then allows us to bring in \cite{HM,HM2} once again. And for groups of Lie type, we consider each possibility for $H_0$ in turn, combining an upper bound on $|{\rm Aut}(H_0)|$ with the well-known dimension bounds of Landazuri and Seitz \cite{LaSe} in order to reduce the problem to a setting where the results of Hiss and Malle are applicable.

The proof of Theorem \ref{t:main} is completed in Section \ref{s:liep}, where we handle the subgroups with socle $H_0 \in {\rm Lie}(p)$. Here we have $H_0 = L/Z(L)$ and 
\[
L = X_{\s} =\{ x \in X \,:\, x^{\s} = x\},
\]
where $X$ is a simple algebraic group over the algebraic closure $k$ of $\mathbb{F}_p$. Then $\bar{V} = V \otimes k$ is an irreducible $kX$-module with highest weight $\l$ and we can appeal to well-known results in the representation theory of simple algebraic groups in order to obtain information on both $n$ and the defining fields for $G_0$ and $H_0$ (for example, see Proposition \ref{p:546_1} in the case where $H_0$ is untwisted). We can then combine this with work of L\"{u}beck \cite{Lubeck}, which gives detailed information on the low-dimensional irreducible $kX$-modules, and this is how we proceed.

\vs

In addition to the appealing simplicity of the bound $|H| < q^n$ in part (iii) of Theorem \ref{t:main}, one of our main motivations for seeking an extension of \cite[Theorem 2.10]{BGS} stems from an application in finite geometry. In order to outline the connection, let $\mathcal{Q} = (\mathcal{P},\mathcal{L})$ be a finite thick generalised quadrangle with point set $\mathcal{P}$ and line set $\mathcal{L}$, and assume $\mathcal{Q}$ has a flag-transitive automorphism group $G$. In this setting, a widely studied conjecture of Kantor \cite{Kantor} from the early 1990s asserts that $\mathcal{Q}$ is either one of the so-called classical generalised quadrangles, or it is one of two sporadic examples, up to duality. In \cite{BFL}, Burness, Feng and Lu complete the proof of Kantor’s conjecture in the special case where $G$ acts primitively on both points and lines.

By earlier work of Bamberg et al. \cite{Bam1}, the latter problem has been reduced to the case where $G$ is an almost simple group of Lie type, which provides  the focus of the recent work in \cite{BFL}. The proof in \cite{BFL} proceeds via a sequence of reductions, which includes the case where $G$ is an almost simple classical group and a point stabiliser $H$ is contained in the collection $\mathcal{S}$. If we write $G_0$ for the socle of $G$ as above, and we set $K = H \cap G_0$, then \cite[Lemma 2.7]{BFL} gives the lower bound
\begin{equation}\label{e:bd4}
|K| > |{\rm Out}(G_0)|^{-5/6}|G_0|^{1/6},
\end{equation}
where ${\rm Out}(G_0)$ is the outer automorphism group of $G_0$. In particular, we can eliminate the candidate pair $(G,H)$ if the lower bound in \eqref{e:bd4} is not satisfied. Of course, in order to implement this approach, we need an upper bound on the order of $H$.

For example, suppose $G_0 = {\rm P\O}_n^{\e}(q)$ is an orthogonal group, where $n \geqs 10$ is even and $q = p^f$ with $p$ a prime. If we set $d = \gcd(4,q^{n/2}-\e)$, then  
\[
|G_0| = \frac{1}{d}q^{\frac{1}{4}n(n-2)}(q^{n/2}-\e)\prod_{i=1}^{(n-2)/2}(q^{2i}-1),\;\; |{\rm Out}(G_0)| = 2df
\]
(see \cite[Table 5.1.A]{KL}) and one can check that  
\[
q^{2n+4} \leqs |{\rm Out}(G_0)|^{-5/6}|G_0|^{1/6}
\]
only if $n \geqs 28$. Therefore, the upper bound $|H|<q^{2n+4}$ in \cite[Theorem 2.10]{BGS} is not sufficient for the desired application in \cite{BFL}. 

However, we find that $q^n \leqs |{\rm Out}(G_0)|^{-5/6}|G_0|^{1/6}$ for all $n \geqs 16$. So putting the special cases arising in parts (i) and (ii) of Theorem \ref{t:main} to one side (they are treated separately in \cite{BFL}), this allows us to assume $n \in \{10,12,14\}$ and we can then consult the relevant tables in \cite[Chapter 8]{BHR} and \cite[Chapter 11]{Schroder} for a complete description of the subgroups in $\mathcal{S}$. In this way, 
Theorem \ref{t:main} plays an important role in the proof of the main theorem of \cite{BFL}, and we anticipate that it will find other applications in future work (for example, in ongoing efforts to extend the results in \cite{BFL} to generalised hexagons and octagons).

\vs

\noindent \textbf{Notation.} Our group-theoretic notation is standard. In particular, we adopt the notation from \cite{KL} for simple groups of Lie type. All logarithms in this paper are in base $2$, unless stated otherwise. Given a pair of positive integers $a,b$, we use the familiar Kronecker delta notation $\delta_{a,b}$, in which case 
$\delta_{a,b} = 1$ if $a=b$, otherwise $\delta_{a,b} = 0$.

\vs

\noindent \textbf{Acknowledgements.} Lei Wang thanks the China Scholarship Council for their financial support, and he thanks the School of Mathematics at the University of Bristol for hosting a research visit in 2026. He also acknowledges the support of NSFC Grants no. 12571022 and 12061083.

\vs

\noindent \textbf{AI statement.} None of the main results or proofs in this paper were obtained through the use of AI tools. In addition, the manuscript has been written without any AI assistance, except for proofreading of the final draft with the aid of ChatGPT 6 Astra.

\section{Classical groups and the collection $\mathcal{S}$}\label{s:prel}

Let $G$ be an almost simple classical group over $\mathbb{F}_{q}$ with socle $G_0$ and natural module $V$. Set $n = \dim V$ and write $q = p^f$, where $p$ is a prime. In view of the existence of exceptional isomorphisms among certain low-dimensional classical groups, as detailed in \cite[Proposition 2.9.1]{KL}, we may assume $G_0$ is one of the following groups:
\begin{equation}\label{e:groups}
\begin{array}{ll}
\mbox{Linear:} \; & \hspace{4.4mm} {\rm L}_n(q), \; \mbox{$n \geqs 2$} \\
\mbox{Unitary:} \; & \hspace{4.4mm} {\rm U}_n(q), \; \mbox{$n \geqs 3$} \\
\mbox{Symplectic:} \; & \hspace{4.4mm} {\rm PSp}_n(q)', \; \mbox{$n \geqs 4$ even} \\
\mbox{Orthogonal:} \; & \left\{ \begin{array}{l}
\O_n(q),\; \mbox{$n \geqs 7$, $nq$ odd} \\
{\rm P\O}_{n}^{\pm}(q),\; \mbox{$n \geqs 8$ even}
\end{array}\right.
\end{array}
\end{equation}
In addition, as explained in Remark \ref{r:main}(a), we will assume that
\begin{equation}\label{e:excep}
G_0 \ne {\rm L}_2(5), \; {\rm L}_2(7), \; {\rm Sp}_4(2)', \; {\rm PSp}_4(3).
\end{equation}

Let  $H$ be a core-free maximal subgroup of $G$. By a theorem of  Aschbacher \cite{asch}, either $H$ belongs to one of eight subgroup collections defined in terms of the geometry of $V$, or $H$ is almost simple and acts irreducibly on $V$. The latter collection is denoted by $\mathcal{S}$ and we will adopt the following formal definition (see \cite[p.3]{KL}). Notice that the various conditions appearing in the definition are designed to ensure that a subgroup in $\mathcal{S}$ is not contained in one of the eight geometric subgroup collections.

\begin{defn}\label{d:sdef}
Suppose $G_0 \ne {\rm P\O}_8^{+}(q)$, ${\rm Sp}_4(q)'$ (with $q$ even) and let $H$ be a core-free maximal subgroup of $G$. Then $H \in \mathcal{S}$ if and only if all of the following conditions are satisfied:
\begin{itemize}\addtolength{\itemsep}{0.2\baselineskip}
\item[(i)] The socle $H_0$ of $H$ is a nonabelian simple group and $H_0 \not\cong G_0$.
\item[(ii)] If $L$ is the full covering group of $H_0$, and if $\rho : L \to {\rm GL}(V)$ is a representation of $L$ such that $\rho(L)=H_0$ (modulo scalars), then $\rho$ is absolutely irreducible.
\item[(iii)] $\rho(L)$ cannot be realised over a proper subfield of $\mathbb{F}$, where 
$\mathbb{F} = \mathbb{F}_{q^2}$ if $G_0={\rm U}_{n}(q)$ is a unitary group, otherwise $\mathbb{F} = \mathbb{F}_{q}$.
\item[(iv)] If $\rho(L)$ fixes a nondegenerate quadratic form on $V$ then $G_0={\rm P\Omega}_{n}^{\e}(q)$.
\item[(v)] If $\rho(L)$ fixes a nondegenerate alternating form on $V$, but no nondegenerate quadratic form, then $G_0={\rm PSp}_{n}(q)$.
\item[(vi)] If $\rho(L)$ fixes a nondegenerate hermitian form on $V$ then $G_0={\rm U}_{n}(q)$.
\item[(vii)] If $\rho(L)$ does not fix a form as in (iv), (v) or (vi) then $G_0={\rm L}_{n}(q)$.
\end{itemize}
\end{defn}

We need to modify the definition of the collection $\mathcal{S}$ in the special cases where $G_0 = {\rm P\O}_8^{+}(q)$ or ${\rm Sp}_4(q)$ (with $q$ even in the latter case). Indeed, both cases were excluded in the statement of Aschbacher's main theorem in \cite{asch} due to complications arising from the existence of special automorphisms. For $G_0 = {\rm Sp}_4(q)$ with $q$ even we refer the reader to \cite[Section 14]{asch} for a description of the maximal subgroups of $G$, while the groups with socle $G_0 = {\rm P\O}_8^{+}(q)$ were handled in later work of Kleidman \cite{Kleidman} (also see \cite[Tables 8.14, 8.50]{BHR}). For the purposes of this paper, we extend the definition of the collection $\mathcal{S}$ as follows:
\begin{itemize}\addtolength{\itemsep}{0.2\baselineskip}
\item[(a)] For $G_0 = {\rm Sp}_4(q)$ with $q  = 2^f$ even, we say that $H \in \mathcal{S}$ if and only if $f \geqs 3$ is odd and $H \cap G_0 = {}^2B_2(q)$.
\item[(b)] Now assume $G_0 = {\rm P\O}_8^{+}(q)$. By inspecting  \cite[Table 8.50]{BHR}, we say that $H \in \mathcal{S}$ if and only if $H \cap G_0$ is one of the following:
\[
\hspace{12mm} \O_7(q) \, (p \geqs 3), \; {\rm Sp}_6(q) \, (p=2), \; {\rm PGL}_3^{\e}(q) \, (q \equiv \e \imod{3}), \; {\rm P\O}_{8}^{-}(q_0) \, (q = q_0^2),
\]
\[
\hspace{12mm} {}^3D_4(q_0) \, (q = q_0^3),\; \O_8^{+}(2) \, (q=p \geqs 3),\; {}^2B_2(8) \, (q =5), \; A_9 \, (q=2), \; A_{10} \, (q=5).
\]
\end{itemize}

Given this definition, it is worth noting that $\mathcal{S}$ may contain subgroups that are not almost simple in these two special cases. For example, if $G = {\rm P\O}_8^{+}(q).\la \tau \ra$, where $q \equiv 1 \imod{3}$ and $\tau$ is a triality graph automorphism, then $\mathcal{S}$ contains a maximal subgroup $3 \times {\rm PGL}_3(q)$. Similarly, if $G = {\rm Sp}_4(q).\la \gamma \ra$, where $q = 2^f$ with $f \geqs 3$ odd and $\gamma$ is an involutory graph automorphism, then $2 \times {}^2B_2(q)$ is contained in $\mathcal{S}$.

For this reason, we will always refer to $H_0$ as the socle of $H \cap G_0$, which is consistent with our notation in Definition \ref{d:sdef}.
 
\section{Proof of Theorem \ref{t:main}: $H_0 \not\in {\rm Lie}(p)$}\label{s:notlie}

Let $G$ be a finite almost simple classical group over $\mathbb{F}_q$ with socle $G_0$ and write $q = p^f$ where $p$ is a prime. Recall that we assume $G_0$ is one of the groups listed in \eqref{e:groups} and we exclude the handful of groups in \eqref{e:excep}. Recall that for each subgroup $H \in \mathcal{S}$, we write $H_0$ for the socle of $H \cap G_0$.

Our goal in this section is to prove Theorem \ref{t:main} for the subgroups $H \in \mathcal{S}$ with $H_0 \not\in {\rm Lie}(p)$. Therefore, we are handling the cases where $H_0$ is an alternating group, a sporadic group, or a simple group of Lie type over a field $\mathbb{F}_t$ with $\gcd(t,p) = 1$. We will prove the following result.

\begin{thm}\label{t:nondef1}
Let $H \in \mathcal{S}$ be a maximal subgroup of $G$ with $H_0 \not\in {\rm Lie}(p)$. Then either $|H| < q^{n}$ or one of the following holds:
\begin{itemize}\addtolength{\itemsep}{0.2\baselineskip}
\item[{\rm (i)}] $H_0 = A_d$, $q=p$ is a prime and $V$ is the fully deleted permutation module for $H_0$.
\item[{\rm (ii)}] $(G_0,H_0)$ is one of the cases in Table \ref{tab:nondef1}.
\end{itemize}
\end{thm}

\begin{prop}\label{p:nd1}
The conclusion of Theorem \ref{t:nondef1} holds if $H_0$ is a sporadic simple group.
\end{prop}

\begin{proof}
First observe that $|H| \leqs |\mathbb{M}| < 2^{180}$, where $\mathbb{M}$ is the Monster sporadic group, so we may assume $n \leqs 179$. 

For $n \leqs 17$ we can inspect the relevant tables of maximal subgroups in $\mathcal{S}$ which are presented in \cite{BHR} (for $n \leqs 12$), \cite{Schroder} (for $n =13,14,15$) and \cite{Rogers} (for $n = 16,17$). The reader can check that all of the relevant cases that arise have been recorded in Table \ref{tab:nondef1}.

Now assume $18 \leqs n \leqs 179$. Here the possibilities for $(n,H_0)$ have been determined by Hiss and Malle \cite{HM,HM2}, with \cite[Table 2]{HM2} as our main reference. And if $q = p^f$ is a prime power such that the prime $p$ satisfies the restrictions given in the third column of \cite[Table 2]{HM2}, then it is easy to check that $|{\rm Aut}(H_0)| < q^n$ with the exception of the cases recorded in Table \ref{tab:spor1}.

\begin{table}
\[
\begin{array}{lll} \hline
H_0 & (n,p) \\ \hline
{\rm J}_3 & (18,2) \\
{\rm HS} & (20,2) \\
{\rm McL} & (22,2) \\
{\rm Co}_3 & (22,3), (22,2) \\
{\rm Co}_2 & (23,3), (22,2) \\
{\rm Co}_1 & (24,5),(24,3),(24,2) \\
{\rm Ru} &  (28,2) \\
{\rm Fi}_{22} &  (27,2) \\ \hline
\end{array}
\]
\caption{A list of candidate subgroups with $H_0$ sporadic}
\label{tab:spor1}
\end{table}

We now consider each possibility listed in Table \ref{tab:spor1} in turn. In each case, let $\chi$ be the corresponding Brauer character of the appropriate covering group of $H_0$ and let $\iota$ be the Frobenius-Schur indicator of $\chi$. In addition, let $\mathcal{I}$ be the set of algebraic irrationalities of the associated character ring, which is the subring of $\mathbb{C}$ generated by the values of $\chi$. We refer the reader to \cite[Section 4.2]{BHR} for a discussion of Brauer character irrationalities and we adopt the notation therein.

First assume $H_0 = {\rm Co}_2$ and $(n,p) = (23,3)$. According to \cite[Table 2]{HM2}, we have $\iota = +$ and $\mathcal{I}$ is empty, which means that the corresponding representation of $H_0$ is defined over the prime field $\mathbb{F}_3$ and the image preserves a nondegenerate quadratic form. Therefore, $G_0 = \O_{23}(3)$ and this case is recorded in Table \ref{tab:nondef1}. Similarly, if $H_0 = {\rm HS}$ and $(n,p)=(20,2)$ then $G_0 = {\rm Sp}_{20}(2)$ and $|H| > 2^{20}$.

Next assume $H_0 = {\rm Co}_3$, $n = 22$ and $p \in \{2,3\}$. By inspecting \cite[Table 2]{HM2} we see that $\iota = -$ (for $p=2$) and $\iota = +$ (for $p=3$), and in both cases $\mathcal{I}$ is empty. Therefore, $G_0 = {\rm Sp}_{22}(2)$ for $p=2$ and $G_0 = {\rm P\Omega}^{\e}_{22}(3)$ for $p=3$. In the latter case, we can use the fact that $|{\rm Co}_3|$ is divisible by $23$ in order to deduce that $\e=+$. Similarly, if $H_0 = {\rm Co}_2$ and $(n,p) = (22,2)$ then $G_0 = \O_{22}^{+}(2)$, using the divisibility of  $|{\rm Co}_2|$ by $23$ in order to determine the sign. And if $H_0= {\rm Ru}$ and $(n,p) = (28,2)$ then $G_0 = \O_{28}^{-}(2)$, noting that $G_0$ is a minus-type group since $|{\rm Ru}|$ is divisible by $29$. All of these claims concerning the sign of $G_0$ can also be verified using \textsf{OSCAR} \cite{OSCAR} (see \cite[p.257]{BNP}, for example).

The case $H_0 = {\rm McL}$ with $(n,p) = (22,2)$ is similar. Here $G_0 = \O_{22}^{\e}(2)$ and we can identify the natural module for $G_0$ with an irreducible module $W$ for $J={\rm Co}_2$, noting that $H_0$ is a maximal subgroup of $J$ and the restriction of $W$ to $H_0$ is irreducible. This means that $H_0 < {\rm Co}_2 < \O_{22}^{+}(2)$ and thus $\e=+$, as recorded in Table \ref{tab:nondef1} (as above, this can also be checked using \textsf{OSCAR}).

Next suppose $H_0 = {\rm Co}_{1}$, $n=24$ and $p \in \{2,3,5\}$, so $G_0 = {\rm P\O}_{24}^{\e}(p)$. For $p \in \{3,5\}$ we can associate $V$ with the irreducible module $\Lambda \otimes_{\mathbb{Z}} \mathbb{F}_p$ for the double cover $2.H_0$, where $\Lambda$ is the Leech lattice. And then from standard properties of $\Lambda$ it follows that $2.H_0$, and hence $H_0$ itself, preserves a plus-type nondegenerate quadratic form on $V$ and thus $\e=+$. Reduction of the Leech lattice action modulo $2$ gives an embedding of ${\rm Co}_1$ in ${\rm O}_{24}^{+}(2)$, so $\e=+$ when $p=2$.

Now assume $H_0 = {\rm J}_3$, $n=18$ and $p=2$. From \cite[Table 2]{HM2} we see that $\iota = \circ$ and $\mathcal{I} = \{ {\rm z}_3, {\rm b}_5\}$. Set $\mu = {\rm b}_5$ and let $\bar{\mu} \in \mathbb{F}_4$ be its $2$-modular reduction. As recorded in \cite[Table 4.2]{BHR}, we have $\mu \in \mathbb{R}$ and the field extension $\mathbb{F}_2(\bar{\mu})$ has degree $2$ and thus $\bar{\mu} \in \mathbb{F}_4 \setminus \mathbb{F}_2$. The latter property means that $\mathbb{F}_4$ is the minimal field of definition for the corresponding representation, which in turn implies that $G_0 = {\rm L}_{18}(4)$ or ${\rm U}_{18}(2)$ since $\iota =\circ$. Finally, by appealing to \cite[Corollary 4.4.2]{BHR}, noting that $\mu \in \mathbb{R}$ and $\bar{\mu} \not\in \mathbb{F}_2$, we see that $G_0 = {\rm L}_{18}(4)$ and we conclude that $|H|<q^n$. 

The final case $H_0 = {\rm Fi}_{22}$ with $(n,p) = (27,2)$ is very similar. Here $\iota = \circ$ and $\mathcal{I} = \{ {\rm z}_3, {\rm b}_{11} \}$. By inspecting \cite[Table 4.2]{BHR} we deduce that $G_0 = {\rm U}_{27}(2)$ via \cite[Corollary 4.4.2]{BHR}, so this case appears in Table \ref{tab:nondef1}.
\end{proof}

\begin{prop}\label{p:nd2}
The conclusion of Theorem \ref{t:nondef1} holds if $H_0$ is an alternating group.
\end{prop}

\begin{proof}
Write $H_0 = A_d$ with $d \geqs 5$ and note that we are free to assume that $V$ is not the fully deleted permutation module for $H_0$. As before, the cases where $n \leqs 17$ can be read off from the relevant tables in \cite{BHR,Rogers,Schroder}, so we will assume $n \geqs 18$.

First assume $d \geqs 24$. By arguing as in the proof of \cite[Theorem 4.1]{L85}, working with results due to James \cite{James83} and Wagner \cite{Wagner77}, we have $n \geqs d(d-5)/4$ and the result follows since 
\[
\log |H| < d \log d < \frac{1}{4}d(d-5).
\]
(Here, and throughout the paper, all logarithms are in base $2$, unless stated otherwise.) 

Now suppose $5 \leqs d \leqs 23$, so $|H| \leqs 23! < 2^{75}$ and the result follows if $n \geqs 75$. Finally, for $18 \leqs n < 75$ we can read off the possibilities for $(n,H_0)$ by inspecting \cite[Table 2]{HM2} and in every case we find that $|H| < q^n$.
\end{proof}

\begin{prop}\label{p:nd3}
The conclusion of Theorem \ref{t:nondef1} holds if $H_0$ is a group of Lie type. 
\end{prop}

\begin{proof}
Let $H_0$ be a simple group of Lie type defined over $\mathbb{F}_t$, where $\gcd(t,p) = 1$. In view of the previous proposition, we may assume 
\begin{equation}\label{e:exc}
H_0 \not\in \{ {\rm L}_2(4), {\rm L}_2(5), {\rm L}_{2}(9), {\rm L}_4(2), {\rm Sp}_{4}(2)'\}
\end{equation}
since each of these groups is isomorphic to an alternating group. As before, let $\chi$ be the corresponding Brauer character, with Frobenius-Schur indicator $\iota$ and associated set $\mathcal{I}$ of irrationalities. 

Our main tool is a well-known theorem of Landazuri and Seitz \cite{LaSe}, which gives lower bounds on the degrees of non-trivial representations of quasisimple groups of Lie type in their non-defining characteristic. More specifically, we work with the bounds presented in \cite[Table 5.3.A]{KL}, which corrects three small errors  in the original paper \cite{LaSe}. In each case, we read off a lower bound $n \geqs f(H_0)$ and we then compare an upper bound on $|{\rm Aut}(H_0)|$ with $q^{f(H_0)}$, where $q \in \{2,3\}$ is minimal such that $\gcd(t,q) = 1$. 

We now divide the proof into a number of separate cases according to the various possibilities for $H_0$. We begin by assuming $H_0$ is a linear group.

\vs

\noindent \emph{Case 1. $H_0 = {\rm L}_d(t)$ with $d \geqs 3$}

\vs

To begin with, let us assume $(d,t) \ne (3,2),(3,4)$. Then $f(H_0) = t^{d-1}-1$ and we have
\begin{equation}\label{e:pslbd}
|{\rm Aut}(H_0)| < 2\log t \cdot t^{d^2-1}.
\end{equation}
It is routine to check that 
\begin{equation}\label{e:pslbd2}
2\log t \cdot t^{d^2-1} \leqs 2^{f(H_0)}
\end{equation}
unless $(d,t) = (3,3)$, $(5,2)$ or $(6,2)$. 

Suppose $(d,t) = (3,3)$, so $|{\rm Aut}(H_0)| = 11232$ and $n \geqs f(H_0)= 8$. Since $|{\rm Aut}(H_0)| < 2^{14}$ and $\gcd(3,q) = 1$,  we may assume $q=2$ and $8 \leqs n \leqs 13$. Then by  inspecting the relevant tables in \cite{BHR,Schroder}, we see that $G_0 = \O_{12}^{-}(2)$ is the only possibility and this case is recorded in Table \ref{tab:nondef1}. Similarly, if $(d,t) = (5,2)$ then we may assume $(n,q) = (15,3)$, which is ruled out by inspecting \cite{Schroder}. And if $(d,t) = (6,2)$ then $n \geqs 31$ and the result follows since $|H| < 3^{31}$.

To complete the argument in Case 1, we may assume $(d,t) = (3,2)$ or $(3,4)$. For $(d,t) = (3,2)$  we have $f(H_0) = 2$ and  $|H| \leqs 336 < 3^6$, so we may assume $2 \leqs n \leqs 5$. Then by inspecting the relevant tables in \cite[Chapter 8]{BHR}, we see that $|H| \geqs q^n$ if and only if $G_0 = {\rm U}_3(q)$ with $q \in \{3,5\}$. 

Finally, suppose $(d,t) = (3,4)$ and note that $f(H_0) = 4$ and $|{\rm Aut}(H_0)| = 241920$. Since $|H|< 3^{12}$, we may assume $4 \leqs n \leqs 11$. Once again, by consulting \cite[Chapter 8]{BHR} we deduce that $|H| < q^n$ unless $G_0 = {\rm L}_4^{\e}(3)$ or ${\rm U}_6(5)$ (note that if $G_0 = {\rm L}_6(7)$ then \cite[Table 8.25]{BHR} indicates that $|H| \leqs 4|H_0| < 7^6$, which explains why this case is not recorded in Table \ref{tab:nondef1}).

\vs

\noindent \emph{Case 2. $H_0 = {\rm L}_2(t)$}

\vs

In view of \eqref{e:exc}, we may assume $t \ne 4,5,9$ and thus  $f(H_0) = (t-1)/d$ with $d = \gcd(2,t-1)$. It is routine to check that
\[
\log t \cdot t(t^2-1) < 2^{f(H_0)}
\]
for all $t \geqs 41$. And by replacing $\log t$ with $\log_rt$ when $t$ is an $r$-power and $r$ is a prime, we reduce to a handful of cases with $t \leqs 29$ and $n \leqs 15$. By inspecting \cite{BHR, Schroder}, the reader can check that all the cases with $|H| \geqs q^n$ have been recorded in Table \ref{tab:nondef1}.
 
\vs

\noindent \emph{Case 3.} $H_0 = {\rm U}_d(t)$ with $d \geqs 3$

\vs

First assume $d$ is odd, so $f(H_0) = t(t^{d-1}-1)/(t+1)$ according to \cite[Table 5.3.A]{KL}. The upper bound in \eqref{e:pslbd} holds and it is easy to check that the inequality in \eqref{e:pslbd2} is satisfied unless $(d,t) = (3,3), (3,4), (3,5), (5,2)$ or $(7,2)$.

Suppose $(d,t) = (3,3)$. Here $|{\rm Aut}(H_0)| = 12096$, $f(H_0) = 6$ and we deduce that $|H| < q^n$ unless $q=2$ and $n \leqs 13$, or $(n,q) = (6,4)$. By inspecting \cite{BHR,Schroder}, we see that $G_0 = {\rm Sp}_6(2)$ is the only case that arises with $|H| \geqs q^n$.

If $(d,t) = (3,4)$ then $f(H_0) = 12$, $q \geqs 3$ and we check that 
$|{\rm Aut}(H_0)| < 3^{12}$. Similarly, if $(d,t) = (3,5)$ then $f(H_0) = 20$ and $|{\rm Aut}(H_0)| < 2^{20}$, and if $(d,t) = (7,2)$ then $f(H_0) = 42$ and the result follows since $q \geqs 3$.

Finally, suppose $(d,t) = (5,2)$. Here $f(H_0) = 10$ and we may assume $q=3$ and $n \leqs 15$, or $q=4$ and $n \leqs 12$. By inspecting \cite{BHR,Schroder}, we see that $G_0 = {\rm PSp}_{10}(3)$ is the only case that arises, and this is listed in Table \ref{tab:nondef1} since $|H| > 3^{10}$.

Now assume $d \geqs 4$ is even. For now, let us exclude the cases $(d,t) = (4,2), (4,3)$, which means that $f(H_0) = (t^d-1)/(t+1)$. Then by arguing as above, we deduce that $|H| < q^n$ unless $(d,t) = (6,2)$, $q = 3$ and $n=20,21$. By inspecting \cite[Table 2]{HM2}, we see that $G_0 = \O_{21}(3)$ is the only possibility. Now $|H_0| < 3^{21}$ and $2|H_0| > 3^{21}$, so we need to know that $H \ne H_0$ in order to include this case in Table \ref{tab:nondef1}. To see this, let $J = {\rm PGU}_6(2) = H_0.3$ and observe that $J$ has a $21$-dimensional self-dual irreducible representation over $\mathbb{F}_3$ (see \cite{GAPCTL}, for example). Therefore, $J \leqs {\rm SO}_{21}(3) = G_0.2$ and we deduce that $J \leqs G_0$. Therefore, $J \leqs H$ and so this case is recorded in Table \ref{tab:nondef1}.

Suppose $(d,t) = (4,2)$. Here $f(H_0) = 4$ and $|{\rm Aut}(H_0)| < 3^{10}$, so we may assume $n \leqs 9$. By inspecting \cite{BHR}, we see that $|H| \geqs q^n$ when $G_0$ is any one of the following groups:
\[
{\rm U}_4(5), \; {\rm L}_4(7), \; {\rm U}_4(11), \; {\rm U}_5(5), \; {\rm L}_5(7).
\]
In addition, if $G = {\rm L}_4(13).2$ (extended by a suitable graph automorphism), then $H = {\rm U}_4(2).2$ is a maximal subgroup in $\mathcal{S}$ with $|H| \geqs 13^4$, so this latter case is also included in Table \ref{tab:nondef1}.

Finally, suppose $(d,t) = (4,3)$, in which case $f(H_0) = 6$. Here we may assume $q=2$ and $n \leqs 24$, together with a handful of additional cases with $n \leqs 12$ and $q \leqs 17$. For $q = 2$ and $n \geqs 18$ we can inspect \cite[Table 2]{HM2} and we deduce that $G_0 = \O_{20}^{\e}(2)$ is the only possibility. Moreover, we can use \textsf{OSCAR} \cite{OSCAR} to check that $\e=+$. Alternatively, by inspecting the $2$-modular character table of $H_0$ we deduce that $C_V(x) = 0$ for every element $x \in H_0$ of order $9$, which means that $\e=+$ since $C_V(y) \ne 0$ for all $y \in \Omega_{20}^{-}(2)$ of order $9$. This case is recorded in Table \ref{tab:nondef1}. 

For $n \leqs 17$, by inspecting \cite{BHR,Rogers, Schroder}, we see that $G_0 = {\rm L}_6(q)$ with $q \in \{7,13\}$ and $G_0 = {\rm U}_6(q)$ with $q \in \{2,5,11\}$ are the only cases that arise (note that $H \cap G_0 = H_0.2$ when $G_0 = {\rm L}_6(13)$, which implies that $|H| \geqs q^n$ in this case). In order to reduce to this list of cases, all of which are recorded in Table \ref{tab:nondef1}, we need to use the fact that $H$ is non-maximal when $G_0 = {\rm U}_{15}(2)$ (see \cite[Table 11.0.20]{Schroder}). In addition, we also need to observe that $|H| \leqs 4|H_0| < 17^6$ if $G_0 = {\rm U}_{6}(17)$, which can be deduced from the information in \cite[Table 8.27]{BHR}.
 
\vs

\noindent \emph{Case 4.} $H_0 = {\rm PSp}_d(t)$ with $d \geqs 4$

\vs

Write $d=2m$ and first assume $t$ is odd, so $f(H_0) = (t^{m}-1)/2$ and
\[
|{\rm Aut}(H_0)| < \log_3t \cdot t^{d(d+1)/2},
\]
which implies that $|H| < q^n$ unless $(d,t) = (4,3), (4,5),(4,7),(6,3)$ or $(8,3)$. 

First assume $(d,t) = (4,3)$. In view of the isomorphism ${\rm PSp}_4(3) \cong {\rm U}_4(2)$ and our previous work in Case 3, we may assume $p=2$. Since $|{\rm Aut}(H_0)| = 51840$ and $f(H_0)=4$, we may assume $n \leqs 12$, or $q=2$ and $n \in \{13,14,15\}$. By inspecting \cite{BHR}, we see that no cases arise when $n \leqs 12$. And by consulting the relevant tables in \cite{Schroder} we deduce that $G_0 = \O_{14}^{-}(2)$ is the only possibility. However, this case is included in Table \ref{tab:def1}, where we view $H$ as a defining characteristic subgroup with socle $H_0 = {\rm U}_4(2)$.

The case $(d,t) = (4,5)$ is very similar. Here $f(H_0) = 12$ and $|{\rm Aut}(H_0)| < 2^{24}$, so we may assume that $12 \leqs n \leqs 23$. By inspecting \cite[Table 2]{HM2}, we see that no cases arise with $14 \leqs n \leqs 23$. And by consulting the relevant tables in \cite{BHR,Schroder} for $n \in \{12,13\}$, we observe that $q \geqs 4$ and thus $|H|<q^n$.

Next assume $(d,t) = (4,7)$. In this case, we have $f(H_0) = 24$ and we may assume $q=2$ and $n \leqs 28$. By inspecting \cite{HM2}, we reduce to the case $n=24$ and we note that $\iota = \circ$ and $\mathcal{I} = \{ {\rm b}_7 \}$. By consulting \cite[Table 4.2]{BHR} and appealing to \cite[Corollary 4.4.2]{BHR}, we deduce that $G_0 = {\rm L}_{24}(2)$ and this case is  listed in Table \ref{tab:nondef1}.

Next suppose $(d,t) = (6,3)$. Here $f(H_0) = 13$ and we may assume that either $q=2$ and $n \leqs 33$, or $q \in \{4,5\}$ and $n \leqs 16$. By inspecting \cite{Schroder}, we deduce that $G_0 = {\rm U}_{14}(5)$ is a genuine exception when $G = G_0.2$ since  $H = H_0.2$ and $|H|>5^{14}$. And by working with \cite{BHR,HM2,Schroder}, it is straightforward to check that $G_0 = {\rm U}_{13}(2)$ and ${\rm U}_{13}(5)$ are the only other cases that arise here, both of which are recorded in Table \ref{tab:nondef1}.

To complete the proof for $t$ odd, let us assume $(d,t) = (8,3)$. Here $f(H_0) = 40$ and we may assume $q=2$ and $n \leqs 56$. By inspecting \cite{HM2}, we see that $n=40$ is the only possibility. Since $\iota = \circ$ and $\mathcal{I} = \{{\rm z}_3\}$, we can use  \cite[Corollary 4.4.2]{BHR} to deduce that $G_0 = {\rm U}_{40}(2)$ and so this case is listed in Table \ref{tab:nondef1}.

Now suppose $t$ is even, which means that $q$ is odd. Recall that we may assume $(d,t) \ne (4,2)$ since ${\rm Sp}_4(2)' \cong A_6$. If $(d,t) \ne (6,2)$ then 
\[
f(H_0) = \frac{1}{2}t^{m-1}(t^{m-1}-1)(t-1) 
\]
and we deduce that 
\[
|H| \leqs |{\rm Aut}(H_0)| < 2\log t \cdot t^{d(d+1)/2} < q^n.
\]
Finally, suppose $(d,t) = (6,2)$, so  $f(H_0) = 7$ and we may assume that $n \leqs 13$ and $q \leqs 7$. By inspecting \cite{BHR,Schroder} we see that $G_0 = \O_7(q)$ with $q \in \{3,5,7\}$ are the only cases that need to be included in Table \ref{tab:nondef1}.

\vs

\noindent \emph{Case 5.} $H_0 = \O_d(t)$ with $d \geqs 7$ odd and $t$ odd

\vs

Write $d = 2m+1$. If $t \geqs 7$ then $f(H_0) = t^{d-3}-1$ and we get
\[
|H| \leqs |{\rm Aut}(H_0)| < \log_3t \cdot t^{d(d-1)/2} < q^n.
\]
Now assume $t \in \{3,5\}$. If $(d,t) \ne (7,3)$, then $f(H_0) = t^{m-1}(t^{m-1}-1)$ and once again we see that $|H|<q^n$. On the other hand, if  $(d,t) = (7,3)$ then $f(H_0) = 27$ and we just need to consider the possibilities for $G_0$ with $q=2$ and $n \leqs 33$. By inspecting \cite{HM2}, we see that $n=27$ is the only option. Here $\iota = \circ$ and $\mathcal{I} = \{{\rm z}_3\}$, which means that $G_0 = {\rm U}_{27}(2)$ and this case is listed in Table \ref{tab:nondef1}.

\vs

\noindent \emph{Case 6.} $H_0 = {\rm P\O}_d^{\e}(t)$ with $d \geqs 8$ even

\vs

Write $d = 2m$ and first assume $\e=+$. If $t \ne 2,3,5$, then 
$f(H_0) = (t^{m-1}-1)(t^{m-2}+1)$ and it is easy to check that
\begin{equation}\label{e:e3}
|H| \leqs |{\rm Aut}(H_0)| < 6\log t \cdot t^{d(d-1)/2} < q^n.
\end{equation}
Similarly, if $t \in \{2,3,5\}$ and $(d,t) \ne (8,2)$ then 
$f(H_0) = t^{m-2}(t^{m-1}-1)$ and once again we see that $|H|<q^n$. Finally suppose $(d,t) = (8,2)$. Here $f(H_0) = 8$ and we may assume $n \leqs 18$ and $q \leqs 13$. By inspecting \cite{BHR,HM2}, we deduce that $G_0 = {\rm P\O}_8^{+}(q)$ with $q \in \{3,5,7,11,13\}$ are the only cases that arise with $|H| \geqs q^n$. (As recorded in Table \ref{tab:nondef1}, if $q = 11$ then $|H| \geqs q^n$ if and only if $G \ne G_0$. Similarly, if $q = 13$ then we need $G = G_0.S_3$ so that $H = H_0.S_3$ and $|H| \geqs q^n$.)

Now assume $\e=-$. Here $f(H_0) = (t^{m-1}+1)(t^{m-2}-1)$ and the above upper bound on $|H|$ in \eqref{e:e3} is sufficient unless $(d,t) = (8,2)$. In the latter case, we have $f(H_0) = 27$ and it is easy to check that $|{\rm Aut}(H_0)| < 3^{27}$.

\vs

\noindent \emph{Case 7.} $H_0$ is an exceptional group of Lie type

\vs

To complete the proof of the proposition, we may assume $H_0$ is an exceptional group. For $H_0 = E_6(t)$ or ${}^2E_6(t)$ we have $f(H_0) = t^9(t^2-1)$ and the bound $|H| < 2\log t\cdot t^{78}$ is sufficient. Similarly, there are no exceptions if $H_0 = E_7(t)$ or $E_8(t)$ since $f(H_0) = t^{a}(t^2-1)$, where $a=15$ if $H_0 = E_7(t)$, and $a=27$ if $H_0 = E_8(t)$.

Next suppose $H_0 = F_4(t)$. For $t$ odd we have $f(H_0) = t^6(t^2-1)$ and the upper bound $|H| < \log_3t \cdot t^{52}$ is sufficient. Similarly, if $t \geqs 4$ is even then $f(H_0) = t^7(t^3-1)(t-1)/2$ and the result follows. The result also follows when $t=2$, since $f(H_0) = 44$.

Now assume $H_0 = G_2(t)'$. If $t \geqs 5$ then $f(H_0) = t(t^2-1)$ and the bound $|H| \leqs 2\log t \cdot t^{14}$ is sufficient. For $t=3$ we have $f(H_0) = 14$ and we may assume $q=2$ and $n \leqs 23$. By inspecting \cite{HM2,Rogers,Schroder} we see that $G_0 = \O_{14}^{-}(2)$ is the only case that arises with $|H| \geqs q^n$.
Next assume $t=4$. Here $f(H_0) = 12$ and we may assume $q=3$ and $n \leqs 18$, or $(n,q) = (12,5)$, and in the usual manner we find that $G_0 = {\rm PSp}_{12}(5)$ is the only case that arises here with $|H| \geqs q^n$. 

Finally, let us assume $H_0 = G_2(2)'$. Here $H_0 \cong {\rm U}_3(3)$ and so in view of our earlier work in Case 3, we may assume $p = 3$. Since $|{\rm Aut}(H_0)| =12096 < 3^9$, we may assume $n \leqs 8$.  By inspecting \cite{BHR}, we see that $G_0 = {\rm U}_6(3)$ is the only case that arises here and this is included in Table \ref{tab:def1}, where we view $H_0$ as a subgroup in ${\rm Lie}(p)$.

Next suppose $H_0 = {}^3D_4(t)$. Here $f(H_0) = t^3(t^2-1)$ and the upper bound $|H| < 3\log t \cdot t^{28}$ is sufficient. Similarly, if $H_0 = {}^2F_4(t)'$ then $f(H_0) = t^4(t-1)\sqrt{t/2}$ and the upper bound $|H| < \log t\cdot t^{26}$ is sufficient if $t \geqs 8$. For $t=2$ we have $f(H_0) = 16$ and $|{\rm Aut}(H_0)| < 3^{16}$, so there are no exceptions. For $H_0 = {}^2G_2(t)'$ we have $f(H_0) = t(t-1)$ and the bound $|H| < \log_3t \cdot t^7$ is sufficient if $t \geqs 27$. For $t=3$ we have $H_0 \cong {\rm L}_2(8)$ and we need to consider the groups with $n \leqs 10$. By inspecting \cite{BHR}, we see that $G_0 = \O_9(q)$ is the only possibility. However, \cite[Table 8.59]{BHR} gives $q \geqs 13$ in this case and thus $|H|<q^n$.

Finally, suppose $H_0 = {}^2B_2(t)$. If $t \geqs 32$ then $f(H_0) = (t-1)\sqrt{t/2}$ and the result follows since $|H| < \log t\cdot t^5$. For $t=8$ we have $f(H_0) = 8$ and we need to consider the case $q=3$ with $n \leqs 10$. By inspecting \cite{BHR}, we see that no examples arise.
\end{proof}

\vs

This completes the proof of Theorem \ref{t:nondef1}.

\section{Proof of Theorem \ref{t:main}: $H_0 \in {\rm Lie}(p)$}\label{s:liep}

In this final section, we complete the proof of Theorem \ref{t:main} by determining the subgroups $H \in \mathcal{S}$ with $H_0 \in {\rm Lie}(p)$ and $|H| \geqs q^n$, where $H_0$ denotes the socle of $H \cap G_0$.

Fix such a subgroup $H$ and let $H_0$ be the socle of $H \cap G_0$, which we are assuming is a simple group of Lie type over $\mathbb{F}_t$, where $t$ is a $p$-power. For the purposes of this preliminary discussion, let us assume that $H_0$ is not one of the twisted groups ${}^2B_2(t)$, ${}^2G_2(t)'$ or ${}^2F_4(t)'$. 

Recall from Definition \ref{d:sdef} that the embedding of $H_0$ in $G_0$ arises from an absolutely irreducible representation $\rho: L \to {\rm GL}(V)$, where $L$ is an appropriate quasisimple group with $L/Z(L) = H_0$. Let $k$ be the algebraic closure of $\mathbb{F}_p$ and let $X$ be a simple algebraic group over $k$ such that 
\[
X_{\s} = \{x \in X \,:\, x^{\s} = x\} = L,
\]
where $\s$ is a suitable Steinberg endomorphism of $X$. Let $r$ be the rank of $X$ and fix a set $\{\l_1, \ldots, \l_r\}$ of fundamental dominant weights for $X$, which are labelled in the usual manner (see \cite{Bou}, for example). Then $\bar{V} = V \otimes k$ is an irreducible $kX$-module with highest weight $\l = \sum_{i}c_i\l_i$, where each $c_i$ is an integer such that $0 \leqs c_i < t$ (see \cite[Section 5.4]{KL}). 

The restriction of $\bar{V}$ to $L$ is irreducible and we will use the notation $\bar{V} = M(\l)$ to signify that the $kL$-module $\bar{V}$ has highest weight $\l$ (and we will sometimes refer to $\l$ as the highest weight of $V$). We say that $\l$ is \emph{$p$-restricted} if $c_i < p$ for all $i$. If $s$ is a positive integer and $\phi \in {\rm Aut}(L)$ is a field automorphism induced by the map $\mu \mapsto \mu^p$ on $k$, then we write $M(\l)^{(p^s)}$ for the $kL$-module afforded by the representation $\rho\phi^s$. By \cite[Proposition 5.4.2(i)]{KL}, we have $M(\l)^{(p^s)} \cong M(p^s\l)$ as $kL$-modules.

Let $M(\l)^*$ denote the dual module and recall that $G_0$ is a symplectic or orthogonal group if and only if $M(\l) \cong M(\l)^*$ is self-dual (see \cite[Lemma 2.10.15(i)]{KL}). By \cite[Proposition 5.4.3]{KL}, $M(\l)$ is self-dual 
if and only if $w_0(\l) = -\l$, where $w_0$ is the longest element in the Weyl group of $X$.

\begin{itemize}\addtolength{\itemsep}{0.2\baselineskip}
\item[{\rm (a)}] Suppose $M(\l)$ is self-dual, $p \geqs 3$ and $n$ is even. Here we can use a well-known criterion of Steinberg (see \cite[Lemma 79]{Steinberg}) to determine whether or not $G_0$ is symplectic, and we refer the reader to \cite[Table 2.1]{Korhonen}.

\item[{\rm (b)}] Next assume $M(\l)$ is self-dual and $p=2$, so $n$ is even. There is no analogue of \cite[Lemma 79]{Steinberg} in this setting and determining whether or not $G_0$ is symplectic is an open problem, in full generality. However, for the pairs $(H_0,\l)$ that  will arise below, we can appeal to Korhonen's work in \cite{Korhonen} to resolve the problem in every case.

\item[{\rm (c)}] If $M(\l)$ is self-dual, $n=2m$ is even and $G_0 = {\rm P\O}_n^{\e}(q)$ is orthogonal, then we still have the problem of determining the sign $\e$ of the associated nondegenerate quadratic form on $V$ and there is no straightforward or uniform way to do this, in general. 

In some cases, it may be possible to identify an $m$-dimensional totally singular subspace of $V$, in which case $\e=+$. More generally, if $q$ is odd, then we can determine $\e$ via \cite[Proposition 2.5.10]{KL}, so long as we can calculate the discriminant of the defining quadratic form on $V$, which means that we need to compute the determinant of the Gram matrix of the associated symmetric bilinear form. This approach is not available when $q$ is even, which leads to additional complications, and we do not attempt to provide a precise description of $\e$ in every case arising in Table \ref{tab:def1}. See Remark \ref{r:sl1}, for example.

\item[{\rm (d)}] Finally, suppose $M(\l)$ is not self-dual, in which case $G_0$ is a linear or unitary group, and $V$ is an $n$-dimensional vector space over $\mathbb{F}_s$ with $s = p^e \in \{q,q^2\}$. Of course, if $e$ is odd then $s = q$ and $G_0 = {\rm L}_n(q)$. However, if $e$ is even, then either $s = q$ and $G_0 = {\rm L}_n(q)$, or $s = q^2$ and $G_0 = {\rm U}_n(q)$. In the latter case, we can appeal to \cite[Lemma 2.10.15(ii)]{KL}, which implies that $G_0 = {\rm U}_n(q)$ if and only if $M(\l)^{(q)} \cong M(\l)^*$.
\end{itemize}

Our main result is the following.

\begin{thm}\label{t:def1}
Let $H \in \mathcal{S}$ be a maximal subgroup of $G$ with $H_0 \in {\rm Lie}(p)$. Then either $|H| < q^{n}$ or $(G_0,H_0)$ is one of the cases in Table \ref{tab:def1}.
\end{thm}

We say that a simple group $H_0 \in {\rm Lie}(p)$ over $\mathbb{F}_t$ is \emph{twisted} if $H_0$ is one of the following:
\begin{equation}\label{e:tw}
{\rm U}_d(t), \; {\rm P\O}_d^{-}(t), \; {}^2E_6(t), \; {}^3D_4(t),\; {}^2F_4(t)',\;  {}^2G_2(t)',\; {}^2B_2(t).
\end{equation}
And in all other cases we say that $H_0$ is \emph{untwisted}. 

We begin the proof of Theorem \ref{t:def1} by handling the untwisted groups. To do this, we will frequently apply the following result, which is an immediate consequence of \cite[Proposition 5.4.6(i)]{KL} (also see \cite[Theorem 5.1.13]{BHR}). In the statement, we assume $G_0$ is a finite simple group over $\mathbb{F}_q$ with natural module $V$, where $\dim V = n$ and $q = p^f$ with $p$ a prime. As before, we write $k$ for the algebraic closure of $\mathbb{F}_p$ and $L$ is a quasisimple group with $L/Z(L) = H_0$.

\begin{prop}\label{p:546_1}
Suppose $H \in \mathcal{S}$ and $H_0 \in {\rm Lie}(p)$ is an untwisted simple group of Lie type defined over $\mathbb{F}_{t}$. Then there exists a positive integer $i$ such that $t = q^i$ and $n = \ell^i$, where $\ell$ is the dimension of an irreducible $kL$-module.
\end{prop}

In order to apply Proposition \ref{p:546_1}, we need to study the dimensions of irreducible $kL$-modules and our main source for this information is \cite{Lubeck}. Here L\"{u}beck gives detailed results on the low-dimensional irreducible $kL$-modules, with the relevant information presented in the form of tables for the low rank groups. This is complemented by \cite[Theorem 5.1]{Lubeck} for the remaining large rank groups, which determines all the irreducible $kL$-modules of dimension at most $r^3$ (up to a small constant), where $r$ is the rank of $X$, as above. It is worth pointing out that L\"{u}beck only works with  $p$-restricted modules in \cite{Lubeck}, but this is sufficient for a complete description of the relevant modules via Steinberg's twisted tensor product theorem (see \cite[Theorem 5.4.5]{KL}).

\begin{prop}\label{p:psldef}
The conclusion of Theorem \ref{t:def1} holds if $H_0 = {\rm L}_d(t)$ and $d \geqs 2$.
\end{prop}

\begin{proof}
Let $H_0 = {\rm L}_d(t)$ and note that $r = d-1$ is the rank of $X$.  As above, let $\l = \sum_i c_i\l_i$ be the highest weight of $V$ and note that $w_0$ acts on the root system of $X$ as $-\tau$, where $\tau$ is the standard graph automorphism. It follows that $V$ is self-dual if and only if $c_{i} = c_{d-i}$ for all $i$.

By Proposition \ref{p:546_1} we have $t = q^i$ and $n = \ell^i$, where $\ell \geqs d$ is the dimension of an irreducible $kL$-module. Therefore,
\[
|H| < i(2-\delta_{2,d})\log q\cdot q^{i(d^2-1)}
\]
and it is easy to check that $|H|<q^{d^i} \leqs q^n$ for all $i \geqs 4$. And similarly if $i=3$ and $d \geqs 3$, and also if $(i,d) = (3,2)$ and $\ell \geqs 3$. If $(i,d,\ell) = (3,2,2)$ then $n = 8$ and by inspecting \cite{BHR} we see that $G_0 = {\rm PSp}_8(q)$ with $q$ odd. This case is recorded in Table \ref{tab:def1}.

Next assume $i=2$, so $H_0 = {\rm L}_d(q^2)$. Suppose $\ell=d$, so $n = d^2$. In view of  \cite[Proposition 5.4.6(i)]{KL}, we may assume $\l = (q+1)\l_1$ or $\l_1+q\l_r$, up to quasiequivalence. In the latter case, we observe that $M(\l)^{(q)} \cong M(\l)^*$ and thus $G_0 = {\rm U}_n(q)$ by \cite[Lemma 2.10.15(ii)]{KL}, whereas $G_0 = {\rm L}_n(q)$ when $\l = (q+1)\l_1$. Both of these cases are recorded in Table \ref{tab:def1}. 

Now assume $i=2$ and $\ell>d$. If $d=2$ then $n \geqs 9$ and $|H|<q^n$. Similarly, if $d=3$ then $n \geqs 36$ (see \cite[Table A.6]{Lubeck}) and the result follows. And if $d \geqs 4$ then $\ell \geqs d(d-1)/2$ (see \cite{Lubeck}) and once again we deduce that $|H|<q^n$.

Finally, suppose $i=1$ and note that $\ell > d$. If $n \geqs d^2$ then
\[
|H| < 2\log_pq \cdot q^{d^2-1} \leqs q^{n}
\]
since $2\log_pq \leqs q$. Therefore, we may assume $n<d^2$. In particular, by applying a field twist, if necessary, we are free to assume that the highest weight $\l$ of $V$ is $p$-restricted. Let $W$ be the natural module for $H_0$.

If $d \geqs 13$ then \cite[Theorem 5.1]{Lubeck} implies that $\l$ is one of the following, up to duals and field twists, and all of these cases are listed in Table \ref{tab:def1}: 
\begin{itemize}\addtolength{\itemsep}{0.2\baselineskip}
\item[{\rm (a)}] $\l = \l_2$ and $n = d(d-1)/2$.
\item[{\rm (b)}] $\l = 2\l_1$, $n = d(d+1)/2$ and $p \geqs 3$.
\item[{\rm (c)}] $\l = \l_1+\l_{d-1}$ and $n = d^2-c_1$,  
where $c_1 = 2$ if $p$ divides $d$, otherwise $c_1=1$. 
\end{itemize}
In (a) and (b), we see that $V$ is not self-dual as a $kL$-module, and by 
appealing to \cite[Lemma 2.10.15(ii)]{KL} we deduce that $G_0 = {\rm L}_n(q)$. In (c), $V$ is self-dual and thus $G_0$ is a symplectic or orthogonal group. We claim that 
\[
G_0 = \left\{\begin{array}{ll} {\rm Sp}_n(q) & 
\mbox{if $p=2$ and $d \equiv 2 \imod{4}$} \\ {\rm P\O}_n^{\e}(q) & \mbox{otherwise.} 
\end{array}\right.
\]
This follows from Steinberg's criterion \cite[Lemma 79]{Steinberg} when $p \ne 2$, and we can appeal to \cite[Theorem 5.1]{Korhonen} when $p=2$. We refer the reader to Remark \ref{r:sl1} for further information on the sign $\e$ when $G_0 = {\rm P\O}_n^{\e}(q)$ with $n$ even. 

Now assume $d \leqs 12$. If $d=2$ then $n \geqs 4$ and the bound $|H|<q^n$ is clear. Now suppose $d \geqs 3$ and note that the possibilities for $n$ are recorded in \cite[Tables A.6-A.15]{Lubeck}. By inspection, we see that $\l = \l_3$ with $d \in \{6,7,8\}$ are the only cases that arise (in addition to the cases labelled (a), (b) and (c) above). Here $n = d(d-1)(d-2)/6$ and we have $G_0 = {\rm L}_n(q)$ if $d = 7,8$, whereas $V = \Lambda^3(W)$ is self-dual when $d=6$. In the latter case, by inspecting \cite[Table 2.1]{Korhonen} we see that $G_0 = {\rm PSp}_{20}(q)$ when $p \ne 2$, whereas \cite[Theorem 4.2]{Korhonen} implies that $G_0 = \O_{20}^{\e}(q)$ when $p=2$. Moreover, in the latter case we can find a $10$-dimensional totally singular subspace of $V = \Lambda^3(W)$ with respect to the defining quadratic form and thus $\e = +$. All of these extra cases are listed in Table \ref{tab:def1}. Note that if $d = 4$ and $\l = \l_2$ then $n = 6$ and this case is excluded because $H_0 \cong G_0 = {\rm P\O}_6^{+}(q)$.
\end{proof}

\begin{rem}\label{r:sl1}
Let us look more closely at the case where $H_0 = {\rm L}_d(q)$ and $\l = \l_1+\l_{d-1}$, which arose in the proof of the previous proposition. We showed that $G_0 = {\rm P\O}_n^{\e}(q)$ when $p \geqs 3$ or $d \equiv 0 \imod{4}$, and here we explain how the sign of $\e$ can be determined when $n$ is even. 

If $\gcd(p,d) = 1$ and $d$ is odd, then $n = d^2-1$ is even, $G_0 = {\rm P\O}_n^{\e}(q)$ and \cite[Lemma 5.4.12]{BHR} implies that $\e=+$ if and only if one of the following holds:
\begin{itemize}\addtolength{\itemsep}{0.2\baselineskip}
\item[{\rm (i)}] $p$ is odd and $(-1)^{(d-1)/2}d$ is a square in $\mathbb{F}_q$; or
\item[{\rm (ii)}] $p=2$ and either $f$ is even or $d\equiv\pm 1\imod{8}$.
\end{itemize}

Now assume $p$ divides $d$ and $d$ is even, so $n = d^2-2$. If $p \ne 2$ then $G_0 = {\rm P\O}_n^{\e}(q)$ and by appealing to \cite[Lemma 5.4.11]{BHR} we may write $V = V_1 \perp V_2$, where $V_1$ and $V_2$ are nondegenerate subspaces with $V_1$ of plus-type and $\dim V_2 = d-2$. The Gram matrix of the associated symmetric bilinear form on $V_2$ has determinant $-1$ and thus $V_2$ is a plus-type space if and only if $(-1)^{d/2}$ is a square in $\mathbb{F}_q$. Therefore, $\e=+$ if and only if $q \equiv 1 \imod{4}$ or $d \equiv 0 \imod{4}$.

Now assume $p=2$ and $d \equiv 0 \imod{4}$, so $G_0 = \O_n^{\e}(q)$. Once again, we have an orthogonal decomposition $V = V_1 \perp V_2$ with $V_1$ plus-type and $\dim V_2 = d-2$. If we write $d = 4m$, then we can decompose $V_2$ as an orthogonal sum of $2m-1$ nondegenerate $2$-spaces, $m-1$ of which are plus-type and $m$ are of type $\delta = (-1)^f$. We conclude that $\e = +$ if and only if $f$ is even or $d \equiv 0 \imod{8}$.
\end{rem}

\begin{prop}\label{p:pspdef}
The conclusion of Theorem \ref{t:def1} holds if $H_0 = {\rm PSp}_d(t)$ and $d \geqs 4$.
\end{prop}

\begin{proof}
Here $w_0$ acts on the root system of $X$ as $-1$, so $w_0(\l) = -\l$ and thus $V$ is self-dual. In particular, $G_0$ is either symplectic or orthogonal.

As before, Proposition \ref{p:546_1} implies that $t = q^i$ and $n = \ell^i$, where $i \geqs 1$ and $\ell \geqs d$ is the dimension of a nontrivial irreducible $kL$-module. Therefore,  
\[
|H| < 2i\log q\cdot q^{id(d+1)/2}
\]
and it is easy to check that $|H|<q^{d^i} \leqs q^n$ for all $i \geqs 3$. 

Next suppose $i=2$. If $\ell > d$ then $\ell \geqs d+2$ (see \cite{Lubeck}) and the bound $|H|<q^n$ holds. Now assume $\ell=d$, so $\l = (1+q)\l_1$ and $n = d^2$. Let $U$ be the natural module for $L = {\rm Sp}_d(q^2)$ and observe that $L$ preserves a nondegenerate quadratic form $Q$ of plus-type on $U \otimes U^{(q)}$, when viewed as a vector space over $\mathbb{F}_{q^2}$ (see \cite[Proposition 1.9.4]{BHR}). If $\e$ denotes the sign of the corresponding $\mathbb{F}_q$-form on $V$, then $G_0 = {\rm P\O}_n^{\e}(q)$ and we calculate that $\e = (-1)^{d/2}$. Note that if $(d,p)=(4,2)$ then $H_0 < {\rm Sp}_8(q) < G_0$ and thus $H$ is non-maximal  (see \cite[Table 7.8]{Rogers}). The cases with $(d,p) \ne (4,2)$ are recorded in Table \ref{tab:def1}.

For the remainder, let us assume $i=1$, so $\ell > d$. Write $d = 2m$. If $n \geqs d^2$ then $|H|<q^n$ and so we may assume $n<d^2$. As before, this allows us to assume that $V$ is $p$-restricted. 

If $d \geqs 24$ then \cite[Theorem 5.1]{Lubeck} implies that one of the following holds: 
\begin{itemize}\addtolength{\itemsep}{0.2\baselineskip}
\item[{\rm (a)}] $\l = \l_2$ and $n = d(d-1)/2-c_2$, where $c_2 = 2$ if $p$ divides $m$, otherwise $c_2=1$,
\item[{\rm (b)}] $\l = 2\l_1$, $p \ne 2$ and $n = d(d+1)/2$.
\end{itemize}
Both of these cases are recorded in Table \ref{tab:def1}. By applying Steinberg's criterion \cite[Lemma 79]{Steinberg} (for $p \ne 2$) and \cite[Theorem 4.2]{Korhonen} (for $p=2$), we deduce that 
$G_0 = {\rm Sp}_n(q)$ if and only if $\l = \l_2$ with $p=2$ and $d \equiv 4 \imod{8}$. If $n$ is even in (b), then $d \equiv 0 \imod{4}$ and we calculate that the discriminant of the defining quadratic form on $V$ is a square in $\mathbb{F}_q$, whence $G_0 = {\rm P\O}_n^{+}(q)$ if and only if $d \equiv 0 \imod{8}$ or $q \equiv 1 \imod{4}$. 

Next assume $6 \leqs d \leqs 22$. By inspecting \cite[Tables A.32-A.40]{Lubeck}, we see that the following cases arise (in addition to cases (a) and (b) given above):
\begin{itemize}\addtolength{\itemsep}{0.2\baselineskip}
\item[{\rm (c)}] $\l = \l_{d/2}$, $d \in \{6,8,10,12\}$, $n = 2^{d/2}$, $p = 2$ and $G_0 = \O_n^{+}(q)$.
\item[{\rm (d)}] $\l = \l_3$, $d = 6$, $n = 14$, $p \geqs 3$ and $G_0 =  {\rm PSp}_{n}(q)$.
\end{itemize}
In case (c), the group $G_0$ is given in \cite[Proposition 5.4.9(i)]{KL}, whereas for (d) we can apply Steinberg's criterion to see that $G_0$ is symplectic. 
Note that if $d=10$ in (c), then $H_0 < \O_{12}^{+}(q) < G_0$ and so $H$ is not maximal. All of the other cases are recorded in Table \ref{tab:def1}.

Finally, by inspecting \cite[Table A.22]{Lubeck}, we see that no further cases arise when $d=4$.
\end{proof}

\begin{prop}\label{p:orth_odddef}
The conclusion of Theorem \ref{t:def1} holds if $H_0 = \O_d(t)$ and $d \geqs 7$ is odd.
\end{prop}

\begin{proof}
Once again, $V$ is self-dual and so $G_0$ is either symplectic or orthogonal. In addition, $t = q^i$ and $n = \ell^i$ as before. Write $d=2m+1$ and note that 
\begin{equation}\label{e:bd6}
|H| < i\log_3q\cdot q^{id(d-1)/2}
\end{equation}
since $q$ is odd. This implies that $|H|<q^{d^i} \leqs q^n$ for all $i \geqs 2$, so we may assume that $i=1$ and $n = \ell>d$.

If $n \geqs d(d-1)/2+1$ then the upper bound in \eqref{e:bd6} yields $|H|<q^n$. Now suppose $d < n \leqs d(d-1)/2$. If $d \geqs 15$ then by inspecting \cite{Lubeck} we see that $\l = \l_2$ is the only possibility. Here $n = d(d-1)/2$ and Steinberg's criterion gives $G_0 = {\rm P\O}_n^{\e}(q)$. By identifying $V$ with the exterior square $\Lambda^2(U)$, where $U$ is the natural module for $H_0$, we calculate that $\e=-$ if and only if $d \equiv 5 \imod{8}$ and $q \equiv 3 \imod{4}$. In addition, a calculation with the spinor norm shows that ${\rm SO}_d(q) = H_0.2$ does not embed in $G_0$, whence $H \cap G_0 = H_0$ and we deduce that $|H| \geqs q^n$ if and only if $|G:G_0|>2$, as recorded in Table \ref{tab:def1}. 

Finally, let us assume $7 \leqs d \leqs 13$ and note that we can read off the possibilities for $n$ from \cite[Tables A.23-A.26]{Lubeck}. In addition to the case $\l = \l_2$ discussed in the previous paragraph, we see that $\l = \l_m$ with $n = 2^m$ is the only other possibility. Here $G_0 = {\rm PSp}_n(q)$ if $d \in \{11,13\}$, otherwise $G_0 = {\rm P\O}_n^{+}(q)$ (see \cite[Proposition 5.4.9(i)]{KL}). For $d = 11$ we note that  $H_0 < {\rm P\O}_{12}^{+}(q) < {\rm PSp}_{32}(q)$ and $H$ is non-maximal. All of the other cases are recorded in Table \ref{tab:def1}.
\end{proof}

\begin{prop}\label{p:orth_plusdef}
The conclusion of Theorem \ref{t:def1} holds if $H_0 = {\rm P\O}_d^{+}(t)$ and $d \geqs 8$ is even.
\end{prop}

\begin{proof}
Write $d = 2m$ and $\l = \sum_i c_i\l_i$. Note that $V$ is self-dual unless $m$ is odd and $c_{m-1} \ne c_m$. As before, we have $t = q^i$ and $n = \ell^i$. We have   
\[
|H| < 2(1+2\delta_{8,d})i\log q\cdot q^{id(d-1)/2},
\]
which implies that $|H|<q^{d^i} \leqs q^n$ for all $i \geqs 2$. 

For the remainder, we may assume $i=1$. If $n \geqs d(d-1)/2+2$ then $|H|< q^n$, so we may assume $n \leqs d(d-1)/2+1$. If $d \geqs 24$ then \cite[Theorem 5.1]{Lubeck} implies that $\l = \l_2$ and $n = d(d-1)/2-c$, where $c = 0$ if $p \geqs 3$ and $c = \gcd(2,m)$ if $p=2$. Here $V$ is self-dual and by combining Steinberg's criterion \cite[Lemma 79]{Steinberg} with \cite[Theorem 4.2]{Korhonen} we deduce that $G_0 = {\rm PSp}_n(q)$ if $d \equiv 4 \imod{8}$ and $p=2$, otherwise $G_0 = {\rm P\O}_n^{\e}(q)$. We only record the case $p \ne 2$ in Table \ref{tab:def1} because $H$ is not maximal when $p=2$ (since we have $H_0 < {\rm Sp}_d(q) < G_0$). In addition, we calculate that $\e = -$ if and only if $d \equiv 4 \imod{8}$ and $q \equiv 3 \imod{4}$.

Finally, suppose $8 \leqs d \leqs 22$. Consulting \cite[Tables A.41-A.48]{Lubeck} we see that the only additional case that arises is $\l = \l_{d/2-1}$ with $d \in \{10,12,14\}$ and $n = 2^{d/2-1}$. Here we have $G_0 = {\rm L}_n(q)$ when $d =10,14$ and 
\[
G_0 = \left\{\begin{array}{ll} 
{\rm PSp}_{n}(q) & \mbox{if $p \geqs 3$} \\
{\rm P\O}_{n}^{+}(q) & \mbox{if $p=2$}
\end{array}\right.
\]
when $d = 12$ (see \cite[Proposition 5.4.9]{KL}).
\end{proof}

\begin{prop}\label{p:exdef1}
The conclusion of Theorem \ref{t:def1} holds if $H_0 = G_2(t)$, $F_4(t)$, $E_6(t)$, $E_7(t)$ or $E_8(t)$. 
\end{prop}

\begin{proof}
As before, by appealing to Proposition \ref{p:546_1} we see that $t = q^i$ and $n = \ell^i$, where $\ell$ is the dimension of a nontrivial irreducible $kL$-module. Note that $V$ is self-dual unless $H_0 = E_6(t)$ and $\l = \sum_i c_i\l_i$ with $(c_1,c_3) \ne (c_6,c_5)$. 

If $H_0 = E_8(t)$ then $|H| < i\log q.q^{248i}$ and $\ell \geqs 248$, which immediately implies that $|H|<q^n$ unless $i=1$ and $n = 248$, in which case $\l = \l_8$ and $G_0$ is orthogonal (see \cite{Korhonen}, for example). We record this case in Table \ref{tab:def1}. 

Next assume $H_0 = E_7(t)$. Here $|H| < i\log q.q^{133i}$, $\ell \geqs 56$ and by inspecting \cite[Table A.52]{Lubeck} we see that $i = 1$ with $\l \in \{ \l_1,\l_7\}$ are the only cases that arise. For $\l = \l_7$ we have $n=56$ and $G_0 = {\rm PSp}_n(q)$ if $p \geqs 3$ and $G_0 = \O_n^{+}(q)$ if $p=2$ (see \cite[Proposition 5.4.18]{KL}). For $\l = \l_1$ we have $n = 133-\delta_{2,p}$ and by inspecting \cite{Korhonen} we see that $G_0 = \O_{n}(q)$ if $p \geqs 3$, otherwise $G_0 = {\rm Sp}_n(q)$. These cases are listed in Table \ref{tab:def1}.

Now suppose $H_0 = E_6(t)$, so $|H| < 2i\log q.q^{78i}$ and $\ell \geqs 27$. In view of \cite[Table A.51]{Lubeck}, we may assume $i=1$ and $\l \in \{\l_1,\l_2\}$. For $\l = \l_1$ we have $n = 27$ and $G_0 = {\rm L}_{27}(q)$ (see \cite[Proposition 5.4.17]{KL}), whereas for $\l = \l_2$ we have $n = 78-\delta_{3,p}$ and $G_0$ is orthogonal by \cite{Korhonen}. Once again, these cases have been included in Table \ref{tab:def1}. 

Next let us turn to the case $H_0 = F_4(t)$. Here $|H| < 2i\log q.q^{52i}$, $\ell \geqs 25$ and we immediately reduce to the case where $i=1$ and $\l \in \{\l_1,\l_4\}$. For $\l = \l_4$ we have $n = 26-\delta_{3,p}$, whereas $n = 26(2-\delta_{2,p})$ when $\l = \l_1$; both cases are recorded in Table \ref{tab:def1} and we note that $G_0$ is orthogonal by \cite{Korhonen}.

Finally, suppose $H_0 = G_2(t)$. Since $|H| < 2i\log q.q^{14i}$ and $\ell \geqs 6$, we may assume $i=1$ and $n \leqs 15$. Then by inspecting \cite[Table A.49]{Lubeck}, either $\l = \l_1$ and $n = 7-\delta_{2,p}$, or $\l = \l_2$ and $n = 14-7\delta_{3,p}$. Here $G_0 = {\rm Sp}_6(q)$ if $\l = \l_1$ and $p=2$, otherwise $G_0 = {\rm P\O}_n^{\e}(q)$. For $\l = \l_2$, by inspecting \cite{Schroder}, we observe that $H$ is non-maximal when $p = 2$, whereas for $p \geqs 5$ we have $\e = +$ if $q \equiv 1,7 \imod{12}$ and $\e = -$ if $q \equiv 5,11 \imod{12}$.
\end{proof}

To complete the proof of Theorem \ref{t:def1}, we may assume $H_0 \in {\rm Lie}(p)$ is one of the twisted groups in \eqref{e:tw}. 

For the cases $H_0 = {\rm U}_d(t)$, ${\rm P\O}_d^{-}(t)$ or ${}^2E_6(t)$ we will work closely with part (ii) of  \cite[Proposition 5.4.6]{KL} (also see \cite[Theorem 5.1.13]{BHR}). Here it will be helpful to recall that $G_0$ is a finite simple classical group over $\mathbb{F}_q$, so the natural module $V$ for $G_0$ is an $n$-dimensional vector space over $\mathbb{F}_{q^u}$, where $u = 2$ if $G_0 = {\rm U}_n(q)$, otherwise $u=1$. Let us also recall that an irreducible $kL$-module $M(\mu)$ is self-dual if and only if
\begin{itemize}\addtolength{\itemsep}{0.2\baselineskip}
\item[{\rm (a)}] $H_0 = {\rm P\O}_d^{-}(t)$ and $d \equiv 0 \imod{4}$; or 
\item[{\rm (b)}] $\mu$ is fixed under the natural action of the standard graph automorphism $\tau$.
\end{itemize}

\begin{prop}\label{p:psudef}
The conclusion of Theorem \ref{t:def1} holds if $H_0 = {\rm U}_d(t)$ and $d \geqs 3$.
\end{prop}

\begin{proof}
By applying \cite[Proposition 5.4.6(ii)]{KL} we see that one of the following holds:
\begin{itemize}\addtolength{\itemsep}{0.2\baselineskip}
\item[{\rm (i)}] $H_0 = {\rm U}_d(q^i)$, $n = \ell^i$ and $G_0 = {\rm PSp}_n(q)$ or ${\rm P\O}_n^{\e}(q)$, where $\ell$ is the dimension of a nontrivial self-dual irreducible $kL$-module.
\item[{\rm (ii)}] $H_0 = {\rm U}_d(q_0^i)$, $q = q_0^2$, $i$ is odd, $n = \ell^i$ and $G_0 = {\rm L}_n(q)$, where $\ell$ is the dimension of a nontrivial irreducible $kL$-module that is not self-dual.
\item[{\rm (iii)}] $H_0 = {\rm U}_d(q^i)$, $i$ is odd, $n = \ell^i$ and $G_0 = {\rm U}_n(q)$, where $\ell$ is the dimension of a nontrivial irreducible $kL$-module that is not self-dual.
\end{itemize}
So in all cases, we have
\[
|H| < 2i\log q \cdot q^{i(d^2-1)}
\]
and it is easy to check that $|H|<q^{d^i} \leqs q^n$ for all $i \geqs 3$. 

Now assume $i = 2$, so we are in case (i) above. Since no $d$-dimensional irreducible $kL$-module is self-dual, it follows that $\ell>d$ and we can now repeat the argument in the proof of Proposition \ref{p:psldef}.

Finally, let us assume $i=1$. As in the proof of Proposition \ref{p:psldef}, we may assume $n<d^2$ and $\l$ is $p$-restricted. If $d \geqs 13$, then  
by inspecting \cite{Lubeck} we see that the following cases arise:
\begin{itemize}\addtolength{\itemsep}{0.2\baselineskip}
\item[{\rm (a)}] $\l = \l_2$, $n = d(d-1)/2$, $H_0 = {\rm U}_d(q)$ and $G_0 = {\rm U}_n(q)$. 
\item[{\rm (b)}] $\l = 2\l_1$, $n = d(d+1)/2$, $p \geqs 3$, $H_0 = {\rm U}_d(q)$ and $G_0 = {\rm U}_n(q)$.
\item[{\rm (c)}] $\l = \l_1+\l_{d-1}$, $n = d^2-c_1$, $H_0 = {\rm U}_d(q)$ and 
\[
G_0 = \left\{\begin{array}{ll} {\rm Sp}_n(q) & 
\mbox{if $p=2$ and $d \equiv 2 \imod{4}$} \\ {\rm P\O}_n^{\e}(q) & \mbox{otherwise,} 
\end{array}\right.
\]
where $c_1 = 2$ if $p$ divides $d$, otherwise $c_1=1$. 
\end{itemize}
These cases are recorded in Table \ref{tab:def1}. 

Finally, if $3 \leqs d \leqs 12$ then by arguing as in the proof of Proposition \ref{p:psldef}, we see that $|H|<q^n$ unless (a), (b) or (c) holds as above (with $d \geqs 5$ in (a) since $G_0 = {\rm P\O}_6^{-}(q) \cong H_0$ when $d=4$), or $H_0 = {\rm U}_d(q)$, $\l = \l_3$ and $d \in \{6,7,8\}$. In the latter case, we have $n = d(d-1)(d-2)/6$. Moreover, if $d \in \{7,8\}$, then $G_0 = {\rm U}_n(q)$, whereas for $d = 6$ we see that $V$ is self-dual and as in the proof of Proposition \ref{p:psldef} we get $G_0 = {\rm PSp}_{20}(q)$ if $p \geqs 3$ and $G_0 = \O_{20}^{\e}(q)$ if $p=2$. And by calculating with the exterior cube $\Lambda^3(U)$, where $U$ is the natural module for $H_0$, we deduce that $\e=+$ when $p=2$.
\end{proof}

\begin{prop}\label{p:orth_minusdef}
The conclusion of Theorem \ref{t:def1} holds if $H_0 = {\rm P\O}_d^{-}(t)$ and $d \geqs 8$ is even.
\end{prop}

\begin{proof}
Write $d = 2m$. According to \cite[Proposition 5.4.6(ii)]{KL}, one of the following holds:
\begin{itemize}\addtolength{\itemsep}{0.2\baselineskip}
\item[{\rm (i)}] $H_0 = {\rm P\O}_d^{-}(q^i)$, $n = \ell^i$ and $G_0 = {\rm PSp}_n(q)$ or ${\rm P\O}_n^{\e}(q)$, where $\ell$ is the dimension of a nontrivial self-dual irreducible $kL$-module (recall that every $kL$-module is self-dual if $m$ is even).
\item[{\rm (ii)}] $H_0 = {\rm P\O}_d^{-}(q_0^i)$, $q = q_0^2$, $i$ is odd, $n = \ell^i$ and 
\[
G_0 = \left\{ \begin{array}{ll}
{\rm L}_n(q) & \mbox{if $m$ is odd} \\
\mbox{${\rm PSp}_n(q)$ or ${\rm P\O}_n^{\e}(q)$} & \mbox{if $m$ is even,}
\end{array}\right.
\]
where $\ell$ is the dimension of a nontrivial irreducible $kL$-module (which is not self-dual if $m$ is odd). 
\item[{\rm (iii)}] $H_0 = {\rm P\O}_d^{-}(q^i)$, $i$ is odd, $n = \ell^i$ and $G_0 = {\rm U}_n(q)$, where $\ell$ is the dimension of a nontrivial irreducible $kL$-module (which is not self-dual if $m$ is odd). 
\end{itemize}
It follows that 
\[
|H| < 2i\log q\cdot q^{id(d-1)/2},
\]
and it is easy to check that $|H|<q^{d^i} \leqs q^n$ for all $i \geqs 2$. 

Now assume $i=1$. Suppose $d \geqs 24$. By arguing as in the proof of Proposition \ref{p:orth_plusdef}, we may assume $\l = \l_2$ and $n = d(d-1)/2-c$, where $c=0$ if $p \geqs 3$ and $c = \gcd(2,m)$ if $p=2$. Here $H_0 = {\rm P\O}_{d}^{-}(q)$ and by inspecting \cite{Korhonen} we see that $G_0 = {\rm PSp}_n(q)$ if $d \equiv 4 \imod{8}$ and $p=2$, otherwise $G_0 = {\rm P\O}_n^{\e}(q)$. For $p = 2$, we have $H_0 < {\rm Sp}_d(q) < G_0$ and so $H$ is not maximal. The case with $p \geqs 3$ is recorded in Table \ref{tab:def1} and closer scrutiny of the embedding reveals that $\e=+$ if and only if $d \equiv 4 \imod{8}$ and $q \equiv 3 \imod{4}$.

Finally, suppose $8 \leqs d \leqs 22$. In addition to the case $\l = \l_2$ discussed above, we may assume that $\l = \l_{d/2-1}$ with $d \in \{8, 10,12,14\}$ and $n = 2^{d/2-1}$. If $d \in \{10,14\}$ then $H_0 = {\rm P\O}_{d}^{-}(q)$ and $G_0 = {\rm U}_n(q)$, otherwise $q = q_0^2$, $H_0 = {\rm P\O}_{d}^{-}(q_0)$ and either $d=8$ and $G_0 = {\rm P\O}_8^{+}(q)$, or $d=12$ and
\[
G_0 = \left\{\begin{array}{ll}
\O_{32}^{+}(q) & \mbox{if $p=2$} \\
{\rm PSp}_{32}(q) & \mbox{if $p \geqs 3$}
\end{array}\right.
\]
(see \cite[Proposition 5.4.9(iii)]{KL}). If $d=12$, then $H_0 < {\rm P\O}_{12}^{+}(q) < G_0$ and $H$ is not maximal. All of the other cases are listed in Table \ref{tab:def1}.
\end{proof}

\begin{prop}\label{p:exdef2}
The conclusion of Theorem \ref{t:def1} holds if $H_0 = {}^2E_6(t)$. 
\end{prop}

\begin{proof}
By appealing to \cite[Proposition 5.4.6(ii)]{KL}, we deduce that one of the following holds:
\begin{itemize}\addtolength{\itemsep}{0.2\baselineskip}
\item[{\rm (i)}] $H_0 = {}^2E_6(q^i)$, $n = \ell^i$ and $G_0 = {\rm PSp}_n(q)$ or ${\rm P\O}_n^{\e}(q)$, where $\ell$ is the dimension of a nontrivial self-dual irreducible $kL$-module.
\item[{\rm (ii)}] $H_0 = {}^2E_6(q_0^i)$, $i$ is odd, $n = \ell^i$ and $G_0 = {\rm L}_n^{\e}(q)$, where 
\[
q= \left\{\begin{array}{ll}
q_0^2 & \mbox{if $\e=+$} \\
q_0 & \mbox{if $\e=-$}
\end{array}\right.
\]
and $\ell$ is the dimension of a nontrivial irreducible $kL$-module that is not self-dual.
\end{itemize}

Since $|H| < 2i\log_2q.q^{78i}$ and $\ell \geqs 27$ (see \cite[Table A.51]{Lubeck}), we see that $|H|<q^n$ if $i \geqs 2$. Now suppose $i=1$. If $\ell > 78$ then $\ell \geqs 324$ and the result follows, so we may assume $V$ is the minimal or adjoint module for $H_0$. If $\l = \l_2$ then $n = 78-\delta_{3,p}$ and $V$ is self-dual, so we have $H_0 = {}^2E_6(q)$ and $G_0$ is orthogonal (see \cite{Korhonen}). And if $\l = \l_1$ then $n = 27$, $H_0 = {}^2E_6(q)$ and $G_0 = {\rm U}_{27}(q)$ (see \cite[Proposition 5.4.17]{KL}). Both cases are recorded in Table \ref{tab:def1}.
\end{proof}

\begin{prop}\label{p:exdef3}
The conclusion of Theorem \ref{t:def1} holds if $H_0 = {}^2B_2(t)$, ${}^2G_2(t)'$, ${}^2F_4(t)'$ or ${}^3D_4(t)$. 
\end{prop}

\begin{proof}
First assume $H_0 = {}^2F_4(t)'$. Here $t = q^i$ and $G_0 = {\rm Sp}_n(q)$ or $\O_{n}^{\e}(q)$, where $p=2$, $i$ is odd and $n \geqs 26^i$ (see \cite[Remark 5.4.7]{KL}). If $i \geqs 3$ then $|H|<i\log q.q^{26i} < q^n$, so we may assume $i=1$. If $n \geqs 27$ then $|H|<q^n$, so we can assume that $n=26$. Here $G_0 = \O_{26}^{\e}(q)$ by \cite[Proposition 5.4.13]{KL} and we have $H_0 < F_4(q) < G_0$, so $H$ is not maximal.

Now suppose $H_0 = {}^2G_2(t)'$. By appealing to \cite[Remark 5.4.7]{KL} we see that $t = q^i$ and $G_0 = {\rm Sp}_n(q)$ or $\O_{n}^{\e}(q)$, where $p=3$, $i$ is odd and $n \geqs 7^i$. Since $|H| < i\log_3q.q^{7i}$, we may assume $i=1$. If $n \geqs 8$ then $|H|<q^n$, so we may assume $n=7$. Here $G_0 = \O_7(q)$, but $H$ is not maximal since $H_0 < G_2(q) < G_0$. 
 
The case $H_0 = {}^2B_2(t)$ is very similar. Here $t = q^i$ and $G_0 = {\rm Sp}_n(q)$ or $\O_{n}^{\e}(q)$, where $p=2$, $i$ is odd and $n \geqs 4^i$. For $i \geqs 3$, it is easy to see that the upper bound $|H| < i\log q.q^{5i}$ is sufficient. Finally, suppose $i=1$. If $n \geqs 6$ then $|H| < q^n$, so we may assume $n \in \{4,5\}$. By inspecting the relevant tables in \cite[Chapter 8]{BHR}, we see that $G_0 = {\rm Sp}_4(q)$ is the only possibility and this case is recorded in Table \ref{tab:def1}.

Finally, let us assume $H_0 = {}^3D_4(t)$. As a consequence of the discussion in \cite[Remark 5.4.7(a)]{KL}, we have $n = \ell^i$ and $G_0 = {\rm PSp}_n(q)$ or ${\rm P\O}_n^{\e}(q)$, where $\ell$ is the dimension of a nontrivial irreducible $kL$-module. In addition, either 
\begin{itemize}\addtolength{\itemsep}{0.2\baselineskip}
\item[{\rm (i)}] $H_0 = {}^3D_4(q^i)$; or 
\item[{\rm (ii)}] $H_0 = {}^3D_4(q_0^i)$ with $q = q_0^3$ and $i \not\equiv 0 \imod{3}$.
\end{itemize}

In both cases, we have $|H|<3i\log q.q^{28i}$ and the bound $n \geqs 8^i$ is sufficient if $i \geqs 2$. Now assume $i=1$. If $n =8$ then $\l = \l_1$, $H_0 = {}^3D_4(q_0)$, $q = q_0^3$ and $G_0 = {\rm P\O}_8^{+}(q)$, so this case is recorded in Table \ref{tab:def1}. If $n \geqs 35$, then $|H|< q^n$, so we can assume $8<n<35$. By inspecting \cite[Table A.41]{Lubeck}, we see that $\l = \l_2$ with $n = 28-2\delta_{2,p}$ is the only possibility. Here $H_0 = {}^3D_4(q)$ and by \cite{Korhonen} we have $G_0 = {\rm P\O}_{n}^{\e}(q)$. If $p=2$ then $H_0 < F_4(q) < G_0$ and $H$ is not maximal, so we only record the case $p \geqs 3$ in Table \ref{tab:def1}.
\end{proof}

\vs

This concludes the proof of Theorem \ref{t:def1}. By combining this with Theorem \ref{t:nondef1}, we deduce that the proof of Theorem \ref{t:main} is complete.

\end{document}